\documentclass[review]{my-elsarticle}

\usepackage{graphicx}
\usepackage{amssymb}
\usepackage{amsfonts}
\usepackage[overload]{empheq}
\usepackage{amsmath}

\journal{Communications in Nonlinear Science and Numerical Simulation}

\newtheorem{thm}{Theorem}[section]

\newtheorem{defin}{Definition}[section]

\usepackage{lineno}
\modulolinenumbers[5]

\begin{document}

\begin{frontmatter}

\title{Direct transcription methods based on fractional integral approximation formulas 
for solving nonlinear fractional optimal control problems}

\author[aut]{Abubakar Bello Salati}
\ead{a.salati@aut.ac.ir, absalati@gmail.com}

\author[aut]{Mostafa Shamsi}
\ead{m\_shamsi@aut.ac.ir}
 
\author[ua]{Delfim F. M. Torres}
\ead{delfim@ua.pt}
 
\cortext[cor1]{This is a preprint of a paper whose final and definite form is with 
\emph{Commun. Nonlinear Sci. Numer. Simul.}, ISSN 1007-5704, available at
[http://dx.doi.org/10.1016/j.cnsns.2018.05.011].
Submitted 25-March-2018; Accepted for publication 16-May-2018.}

\address[aut]{Department of Applied Mathematics, Faculty of Mathematics and Computer Science,\\ 
Amirkabir University of Technology, No. 424 Hafez Avenue, Tehran, Iran\fnref{label3}}

\address[ua]{Department of Mathematics, Center for Research and Development 
in Mathematics and Applications (CIDMA),\\ 
University of Aveiro, 3810-193 Aveiro, Portugal}


\begin{abstract}
This paper presents three direct methods based on Gr\"{u}nwald-Letnikov, 
trapezoidal and Simpson fractional integral formulas to solve  
fractional optimal control problems (FOCPs). 
At first, the fractional integral form of FOCP is considered, then 
the fractional integral is approximated by  Gr\"{u}nwald-Letnikov, 
trapezoidal and Simpson formulas in a matrix approach. 
Thereafter, the performance index is approximated either 
by trapezoidal or Simpson quadrature. As a result, FOCP are reduced 
to nonlinear programming problems, which can be solved by many well-developed algorithms.
To improve the efficiency of the presented method, the gradient of the
objective function and the Jacobian of constraints are prepared in closed forms.
It is pointed out that the implementation of the methods is simple and, 
due to the fact that there is no need to derive necessary conditions, 
the methods can be simply and quickly used to solve a wide class of FOCPs. 
The efficiency and reliability of the presented methods are assessed by ample 
numerical tests involving  a free final time with path constraint FOCP,
a bang-bang FOCP and an optimal control of a fractional-order HIV-immune system. 
\end{abstract}


\begin{keyword}
Fractional optimal control\sep 
Direct numerical solution\sep  
Fractional integration matrix \sep 
Gr\"{u}nwald-Letnikov, trapezoidal and Simpson fractional integral formulas.
\MSC[2010] 26A33 \sep 49M25.
\end{keyword}

\end{frontmatter}


\section{Introduction}

Fractional calculus may be considered an old and yet an interesting topic \cite{NewMACHADO20111140}. 
It deals with the investigation of integrals and derivatives of an arbitrary order. 
At the initial stage, the fractional derivative was a mathematical tool without 
a tangible application. But presently, fractional calculus has wide applications 
in various disciplines like viscoelastic materials \cite{MR2676137,NewAHMADIAN201744}, 
bioengineering applications \cite{magin2006fractional,NewIONESCU2017141},  
signal processing \cite{marks1981differintegral,MR2768178}, 
mathematical finance \cite{MR1882830,scalas2000fractional}, 
control theory \cite{MR2721980,podlubny1994fractional}
and many other areas \cite{AGARWAL:etal,NewTENREIROMACHADO2018174}.
The nonlocal nature of the fractional calculus has given it a unique characteristic 
in modeling some complex systems with memory and hereditary properties that arise from physical, 
engineering and even economic processes \cite{NewMORELES201781,NewTARASOVA2018127}.  
More elaborate details on the theory and applications of fractional calculus  
can be found in \cite{MR2676137,MR1658022,MR2680847,MR2984893}.

Application of fractional calculus to optimal control problems has given birth 
to another new specialization, known as fractional optimal control 
\cite{MR2984893,MR3443073,MR3331286,MR3714431}.  
According to \cite{MR2595937}, a Fractional Optimal Control Problem (FOCP) 
is an optimal control problem in which the performance index and/or 
the dynamic equations contain at least one fractional 
derivative operator. In simple terms, FOCPs are a generalization of the integer-order 
optimal control problems, which are obtained by replacing integer-order derivatives with fractional ones. 
Recently, some interesting and real-life models of FOCPs have been presented by the researchers in \cite{ding2012optimal,sweilam2016optimal,sweilam2017legendre,NewZAKY2017177,MyID:397}.

Like integer-order optimal control problems, the numerical methods for FOCPs 
can be categorized into  ``direct" and ``indirect" methods.
With indirect methods, the solution is obtained by solving a fractional 
Hamiltonian boundary-value problem, which is derived from the optimality conditions. 
Accordingly, the first step in indirect methods is the derivation 
of first-order optimality conditions.
On the other hand, direct methods do not rely on optimality conditions. 
In these approaches, FOCPs are solved by transcribing them into 
Nonlinear Programming problems (NLP). Thereafter, a NLP-solver 
is used to solve the resulting problem \cite{MR3089375}.
The indirect methods, in comparison with the direct, have some disadvantages, including,
(i) difficulties in deriving the Hamiltonian boundary-value problem, especially 
for problems with path constraints, and
(ii) sensitivity to initial guess for state functions as well as costate or adjoint functions.
Due to the aforementioned reasons, direct methods are easier to employ 
than the indirect methods to solve both integer and fractional order optimal control problems. 


As earlier works on the formulation, derivation of optimality conditions  
and  direct/indirect solution schemes for FOCPs, we can refer to 
\cite{MR2112177,MR2356715,MR2463065,MR2433010,MR2386201,MR2519151,MR2762638}.
After these pioneer works, many other extensive researches have been 
done on the development of numerical methods for FOCPs. For instance, we can refer to 
Oustaloup recursive approximation \cite{MR2595937},
direct methods based on pseudo-state-space formulations of FOCP \cite{MR2849617},
spectral methods based on orthogonal polynomials and fractional operational matrices
\cite{MR2824693,MR3044581,bhrawy2015new,bhrawy2015efficient,doha2015efficient,ezz2017numerical,ezz2018numerical},
Legendre multiwavelet collocation methods \cite{MR2895863},
direct methods based on Bernstein polynomials \cite{MR3179222,keshavarz2015numerical,rabiei2018numerical},
nonstandard finite difference methods \cite{zahra2015non},
linear programming approaches \cite{rakhshan2016efficient},
integral fractional pseudospectral methods \cite{tang2015integral},
direct methods based on Ritz's techniques \cite{jahanshahi2016simple,nemati2016efficient},
the epsilon-Ritz method \cite{JOTALOTFI}, 
direct methods based on hybrid block-pulse with other basis functions
\cite{mashayekhi2016approximate,yonthanthumapproximate},
pseudospectral methods based on Legendre M\"{u}ntz basis functions \cite{ejlali2017pseudospectral},
dynamic Hamilton-Jacobi-Bellman methods \cite{rakhshan2016solving},
penalty and variational methods \cite{lotfi2017combination},
control parameterization methods \cite{mu2017control}, 
differential and integral fractional pseudospectral methods \cite{tang2017new},
as well as other numerical techniques \cite{singha2017efficient,Almeida2014,baleanu2017new}.
Efforts were also done to derive optimality conditions for special types of FOCPs, 
such as bang-bang FOCPs \cite{kamocki2015fractional} and free final and terminal 
time problems \cite{MR3151796,MR3124694}.

Direct transcription methods based on local methods, such as trapezoidal and Simpson methods, 
are widely used to solve integer-order optimal control problems and some softwares, 
such as SOCS \cite{betts1997sparse} and ICLOCS \cite{ICLOCS}, are developed based on these methods. 
However, it is somewhat surprising that these methods have not been yet fully explored to solve FOCPs. 
In this paper, we extend these methods to fractional order optimal control problems.
For this purpose, the matrix form of  Gr\"{u}nwald-Letnikov, trapezoidal and Simpson approximation 
formulas for fractional integrals, which are called \emph{fractional integral matrices}, are derived. 
These fractional integration matrices are used to reduce the FOCP to a NLP. Thereafter, 
the resulted NLP is solved by using a well-developed solver to obtain an approximate solution of the FOCP.
Moreover, in order to increase the speed of the proposed direct methods, the exact gradient 
of the objective function and the Jacobian of constraints are derived and supplied to the NLP-solver.
Finally, the reliability, efficiency and accuracy of the proposed direct method are demonstrated 
with four test problems, including a nonlinear and complex FOCP, a FOCP  with path and terminal 
constraints, a bang-bang FOCP and an applied optimal control problem
of a fractional order HIV-immune system with memory.

The paper is organized as follows. 
In Section~\ref{sec:2}, the notions of fractional derivative and integral are 
reviewed and three approximation methods for the fractional integral, in matrix form, are introduced.
The considered formulation of fractional optimal control problems is stated in Section~\ref{sec:3}. 
The detailed implementation of direct methods is presented in Section~\ref{sec:4}.
In Section~\ref{sec:5}, four numerical examples are provided to show 
the efficiency and reliability of the proposed methods.
Finally, a conclusion is given in Section~\ref{sec:6}.


\section{Definition and approximation of fractional  integral and derivative}
\label{sec:2}

In this section, some definitions, together with three approximation formulas 
for fractional integrals, are reviewed. Moreover, for easy implementation, 
the matrix forms of these approximation formulas are derived.

\begin{defin}[See, e.g., \cite{MR1658022,MR2680847}]
The left sided fractional integral with order $\alpha > 0$ of a 
given function $y(x)$, $x \in (a,b)$, is defined as 
\begin{equation} 
\label{fracintegral}
{}_0\mathcal I^\alpha_x  y\left( x \right) 
= \frac{1}{{\Gamma \left( {\alpha } \right)}} 
\int_0^x {\mathop {\left( {x - t} \right)}\nolimits^{\alpha  - 1} } \mathop y\nolimits \left( t \right)\text{\rm d}t,
\end{equation}
where $\Gamma(\cdot)$ is Euler's gamma function.
\end{defin}

\begin{defin}[See, e.g., \cite{MR1658022,MR2680847}]
The left sided Caputo fractional derivative of order $\alpha\in (0,1)$ 
of a function $y\left( x \right)$ is defined as
\begin{equation}
{}_0^{\text{\tiny C}\hspace{1pt}}\mathcal D^\alpha_x  y\left( x \right) 
= \frac{1}{{\Gamma \left( {1 - \alpha } \right)}}
\int_0^x {\mathop {\left( {x - t} \right)}\nolimits^{-\alpha} } y \left( t \right)\text{\rm d}t.
\end{equation}
\end{defin}

The first step in developing our numerical methods is to approximate the fractional 
integral of a given function. For this purpose, we  briefly review 
Gr\"{u}nwald-Letnikov (GL), trapezoidal (TR) and Simpson (SI) formulas 
for approximating the fractional integral.


\subsection{The Gr\"{u}nwald-Letnikov formula for approximating a fractional integral}

Let the fractional integral of function $y$ be defined on the interval $[0, 1]$. 
To approximate the fractional integral of function $y$, the interval $[0, 1]$ 
is divided into $n$ equal parts by the following mesh points:
\begin{equation*}
x_i:= ih,  \qquad i=0,1,\dots,n,
\end{equation*}
where $h=1/n$. 
The Gr\"{u}nwald-Letnikov (GL) formula, for approximating the fractional integral 
of function $y$ at $x=x_i$, can be expressed as 
\begin{equation} 
\label{GLformula}
\left. {}_0\mathcal I^\alpha_x  y\left( x \right)\right|_{x=x_i} 
\simeq \sum_{j=0}^{i} {}^{\text{\tiny G}}\hspace{-1pt} 
\varpi_{ij}^{(\alpha)}y(x_{j}), \quad i=0,\dots,n,
\end{equation}
where ${}^{\text{\tiny G}}\hspace{-1pt}\varpi_{ij}^{(\alpha)}$ 
are the coefficients of the GL approximation formula, defined by
\begin{equation}
\label{WoGL}
{}^{\text{\tiny G}}\hspace{-1pt}\varpi_{ij}^{(\alpha)}
:=
\begin{cases}
0,& i=0,\\
w_{i-j}^{(\alpha)},& i>0,
\end{cases}
\end{equation}
such that
\begin{equation*}
\omega_{k}^{(\alpha)} :=(-1)^{k}\binom{-\alpha}{k}{h^{\alpha}}
=\frac{(-1)^{k}\Gamma(1-\alpha)}{\Gamma(1+k)\Gamma(1-\alpha-k)}{h^{\alpha}}
\end{equation*}
(see \cite{li2015numerical,NewBRZEZINSKI2016151}). Let 
\begin{subequations}
\label{vectyy}
\begin{align}
&\mathbf y:=\left[y(x_0),y(x_1),\dots,y(x_n)\right]^T,\\ 
&\mathbf{\hat y}^{(-\alpha)} :=\left[{}_0\mathcal I^\alpha_x y(x_0),
{}_0\mathcal I^\alpha_x y(x_1),\dots,{}_0\mathcal I^\alpha_x y(x_n)\right]^T.
\end{align}
\end{subequations}
Then, all the $n+1$ formulas in \eqref{GLformula} can be written 
simultaneously in the following matrix form:
\begin{equation}
\mathbf{\hat y}^{(-\alpha)} \simeq \mathbf W_{\text{\tiny GL}}^{(\alpha)} \,\mathbf y,
\end{equation}
where
\begin{equation}\label{GLmat}
\mathbf W_{\text{\tiny GL}}^{(\alpha)}:=
\begin{bmatrix}
0 & 0 &0&\dots & 0\\
\omega_{1}^{(\alpha)} & \omega_{0}^{(\alpha)} &0& \cdots & 0\\
\omega_{2}^{(\alpha)} &\omega_{1}^{(\alpha)} & \omega_{0}^{(\alpha)} & \cdots & 0\\
\ddots &\ddots & \ddots & \ddots & \vdots\\
\omega_{n}^{(\alpha)} &\omega_{n-1}^{(\alpha)} & \omega_{n-2}^{(\alpha)} & \cdots &  \omega_{0}^{(\alpha)}\\
\end{bmatrix}.
\end{equation}
The matrix  $\mathbf W_{\text{\tiny GL}}^{(\alpha)}$ is called the 
\emph{fractional GL integration matrix}.


\subsection{The trapezoidal formula for approximating a fractional integral}

Let $h=1/n$ and $x_i:=i\hspace{1pt} h,\, i=0,\dots,n$.
If the given function $y$ is approximated  by its piecewise linear 
interpolant based on the nodes $x_i$, $i=0,\dots,n$, and the fractional 
integral is applied to this approximation, then
the trapezoidal (TR) formula  is derived as 
\begin{equation}
\label{TRformula}
\left.{}_0\mathcal I^\alpha_x  y\left( x \right)\right|_{x = x_i} 
\simeq \sum_{j = 0}^{i} {}^{\text{\tiny T}}\hspace{-1pt} 
\varpi_{ij}^{(\alpha)}  y(x_j), \quad i=0,1,\dots,n,
\end{equation}
where
\begin{equation}
\label{WoTR}
{}^{\text{\tiny T}}\hspace{-1pt} \varpi_{ij}^{(\alpha)}:= \frac{h^{\alpha}}{\Gamma (\alpha + 2)}
\begin{cases}
(i - 1)^{\alpha + 1}- (i-1- \alpha)i^{\alpha}, & j=0, i>0, \\
(i - j +1)^{\alpha +1} + (i - 1 - j)^{\alpha + 1} - 2(i -j)^{\alpha +1}, & 1\leq j \leq i-1,\\
1, & j=i,\\
0,&\text{otherwise}
\end{cases}
\end{equation}
(see \cite{li2015numerical}). If we set 
\begin{align*}
& \bar a^{(\alpha)}_k:=
\tfrac{h^{\alpha}}{\Gamma (\alpha + 2)}
\begin{cases}
0,& k=0,\\
(k - 1)^{\alpha + 1}- (k-1- \alpha)k^{\alpha},& k>0,\\
\end{cases}\\
&\bar b^{(\alpha)}_k:=
\tfrac{h^{\alpha}}{\Gamma (\alpha + 2)}
\begin{cases}
1,& k=0,\\
(k +1)^{\alpha +1} + (k - 1)^{\alpha + 1} - 2k^{\alpha +1}, & k>0,
\end{cases}
\end{align*}
then we can express the coefficients of the trapezoidal formula \eqref{TRformula} as
\begin{equation}
{}^{\text{\tiny T}}\hspace{-1pt} \varpi_{ij}^{(\alpha)}:=
\begin{cases}
\bar a^{(\alpha)}_i, & k=0, j>0, \\
\bar b^{(\alpha)}_{i-j}, & 1\leq j \leq i,\\
0,&\text{otherwise}.
\end{cases}
\end{equation}
Considering \eqref{vectyy}, the trapezoidal formula \eqref{TRformula} 
can be expressed in the following matrix form:
\begin{equation}
\mathbf{\hat y}^{(-\alpha)} \simeq \mathbf W_{\text{\tiny TR}}^{(\alpha)} \,\mathbf y,
\end{equation}
where $\mathbf W_{\text{\tiny TR}}^{(\alpha)}$ is the 
\emph{fractional TR integration matrix} defined by
\begin{equation}
\label{TRmat}
\mathbf W_{\text{\tiny TR}}^{(\alpha)}:=
\begin{bmatrix}
0& 0 &0&0&0&\dots &0& 0\\
\bar a^{(\alpha)}_1& \bar b_0^{(\alpha)} &0& 0&0&\cdots & 0&0\\
\bar a^{(\alpha)}_2 &\bar b_1^{(\alpha)} & \bar b_0^{(\alpha)} 
&0& 0&\cdots & 0&0\\
\bar a^{(\alpha)}_3 &\bar b_2^{(\alpha)} & \bar b_1^{(\alpha)} 
& \bar b_0^{(\alpha)}&0&\cdots & 0&0\\
\bar a^{(\alpha)}_3 &\bar b_3^{(\alpha)} & \bar b_2^{(\alpha)} 
& \bar b_1^{(\alpha)}
&\bar b_0^{(\alpha)}&\cdots & 0&0\\
\vdots &\vdots & \vdots & \ddots & \ddots&\ddots&\vdots&\vdots\\
\bar a^{(\alpha)}_n &\bar b_{n-2}^{(\alpha)} & \bar b_{n-3}^{(\alpha)} &\bar b_{n-4}^{(\alpha)} 
&\cdots&\bar b_1^{(\alpha)}  &  \bar b_0^{(\alpha)}&0\\
\bar a^{(\alpha)}_n &\bar b_{n-1}^{(\alpha)} & \bar b_{n-2}^{(\alpha)} &\bar b_{n-3}^{(\alpha)} 
&\cdots&\bar b_2^{(\alpha)} &\bar b_1^{(\alpha)} &  \bar b_0^{(\alpha)}\\
\end{bmatrix}.
\end{equation}


\subsection{The Simpson formula for approximating a fractional integral}

Let $n$ be even and the interval $[0,1]$ be divided into $n$ parts by 
the mesh points $x_i:=i h,\,i=0,\dots,n$.
Then the given function $y$ on the intervals $[x_{2i},x_{2i+2}]$, $i=0,\dots,\tfrac n 2 -2$, 
can be approximated by its corresponding quadratic interpolation polynomial, i.e., 
$y$ is approximated by a piecewise quadratic polynomial on the whole interval $[0,1]$.  
By replacing $y(t)$ with this piecewise quadratic polynomial approximation in \eqref{fracintegral}, 
the Simpson formula for approximating the fractional integral of function $y$ 
on the nodes $x_i$, $i=0,\dots,n$, is derived as
\begin{equation}
\label{SIformula}
\left.{}_0\mathcal I^\alpha_x  y\left( x \right)\right|_{x = x_i} \simeq \sum_{j = 0}^{i+1} 
 {}^{\text{\tiny S}}\hspace{-1pt}\varpi_{ij}^{(\alpha)}
y(x_j), \quad i=0,1,\dots,n,
\end{equation} 
where
\begin{equation}
\label{WoSI}
{}^{\text{\tiny S}}\hspace{-1pt}\varpi_{ij}^{(\alpha)}:=
\begin{cases}
\gamma_{i}^{(\alpha)}, & j=0,\\
\theta_{i-j+1}^{(\alpha)}, &j \text{ is odd},\\	
\mu_{i-j+2}^{(\alpha)}, &j \text{ is even},
\end{cases}
\end{equation}
such that the parameters $\gamma_{k}^{(\alpha)}$, 
$\theta_{k}^{(\alpha)}$ and $\mu_{k}^{(\alpha)}$ are defined as
\begin{align*}
&\gamma_k^{(\alpha)}:=
\frac {{h}^{\alpha}}{\Gamma  \left( \alpha+3 \right) }
\begin{cases}
0,&k\leq 0,\\
\tfrac 1 2\left( 2\,\alpha+3 \right) \alpha,&k=1,\\
\lambda_{0,k}^{(\alpha)}, &2\leq k,
\end{cases}\\
&\theta_k^{(\alpha)}:=
\frac {{h}^{\alpha}}{\Gamma  \left( \alpha+3 \right) }
\begin{cases}
0,&k\leq 0,\\
2\,\alpha+2,&k=1,\\
\lambda^{(\alpha)}_{{1,k}}, &2\leq k,
\end{cases}	\\
&\mu_k^{(\alpha)}:=
\frac {{h}^{\alpha}}{\Gamma  \left( \alpha+3 \right) }
\begin{cases}
0,&k\leq 0,\\
-\tfrac 1 2\alpha,&k=1,\\
\lambda^{(\alpha)}_{{2,2}}, &k=2,\\
\lambda^{(\alpha)}_{{2,3}} +\tfrac 1 2 \left( 2\,\alpha+3 \right) \alpha,&k=3,\\
\lambda^{(\alpha)}_{{2,k}}  +\lambda^{(\alpha)}_{{0,k-2}}, &4\leq k,
\end{cases}    
\end{align*}
where
\begin{align*}
&\lambda_{0,k}^{(\alpha)}:= \frac 1 2\left[ 2\,{\alpha}^{2}- \left( 3\,k-6 \right) 
\alpha+2\,{k}^{2}-6\,k+4 \right] {k}^{\alpha}-\frac 1 2 \left( 2\,k+\alpha-2 \right) 
\left( k-2 \right) ^{\alpha+1},\\
&\lambda_{1,k}^{(\alpha)}:= 2\, \left( k-2 \right) ^{\alpha+1} \left( k+\alpha \right) 
-2\,{k}^{\alpha+1} \left( k-\alpha-2 \right),\\
&\lambda_{2,k}^{(\alpha)}:=\frac 1 2{k}^{\alpha+1} \left( 2\,k-\alpha-2 \right) 
-\frac 1 2\left( k-2\right) ^{\alpha} \left( 2\,{k}^{2}+ \left( 3\,\alpha-2 \right) 
k+2\,{\alpha}^{2} \right).
\end{align*}
Considering \eqref{vectyy}, the Simpson formula \eqref{SIformula} 
can be expressed in the following matrix form
\begin{equation}
\mathbf{\hat y}^{(-\alpha)} \simeq \mathbf W_{\text{\tiny SI}}^{(\alpha)} \,\mathbf y,
\end{equation}
where $\mathbf W_{\text{\tiny SI}}^{(\alpha)}$ is the \emph{fractional SI integration matrix}
\begin{equation}
\label{SImat}
\mathbf W_{\text{\tiny SI}}^{(\alpha)}:=
\begin{bmatrix}
0& 0 &0&0&0&\dots &0& 0\\
\gamma_1^{(\alpha)} & \theta_1^{(\alpha)}     & \mu_1^{(\alpha)}        &     0               
&0                                   &\cdots & 0                              &0\\
\gamma_2^{(\alpha)}       & \theta_2^{(\alpha)}     & \mu_2^{(\alpha)}
&  0   & 0                                  &\cdots & 0                              &0\\
\gamma_3^{(\alpha)}       & \theta_3^{(\alpha)}     & \mu_3^{(\alpha)}        
&  \theta_1^{(\alpha)}     & \mu_1^{(\alpha)} 
&\cdots & 0                               &0\\
\gamma_4^{(\alpha)}       & \theta_4^{(\alpha)}     & \mu_4^{(\alpha)} 
&  \theta_2^{(\alpha)}     & \mu_2^{(\alpha)}     
&\cdots &           0                     &0\\
\vdots  & \vdots  & \vdots  &
\ddots  & \ddots                        &\ddots&    \vdots                   &\vdots\\
\gamma_{n-1}^{(\alpha)} &\theta_{n-1}^{(\alpha)} & \mu_{n-1}^{(\alpha)} & \theta_{n-3}^{(\alpha)}
& \mu_{n-3}^{(\alpha)}&\cdots&\theta_1^{(\alpha)} & \mu_1^{(\alpha)}\\
\gamma_n^{(\alpha)}       &\theta_n^{(\alpha)}      & \mu_n^{(\alpha)}        
& \theta_{n-2}^{(\alpha)}& \mu_{n-2}^{(\alpha)}&\cdots& \theta_2^{(\alpha)}& \mu_2^{(\alpha)}\\
\end{bmatrix}.
\end{equation}
The Simpson and trapezoidal approximation formulas have been presented in \cite{li2015numerical}. 
Due to the fact that these matrix forms have some advantages in the implementation of our method, 
in this paper we derived the explicit forms of these matrices.


\section{The fractional optimal control problem}
\label{sec:3}

In this work, we consider a general formulation for FOCPs, 
which is described as follows: find the optimal control 
$\mathbf u(t)=\left[u_1(t),\dots,u_q(t)\right]\in \mathbb{R}^q$,
the state $\mathbf x(t)=\left[x_1(t),\dots,x_p(t) \right]\in \mathbb{R}^p$ 
and  possibly the terminal time $t_\text{f}\hspace{1pt}$ 
that minimize the performance index
\begin{subequations}
\label{main-focp}
\begin{equation} 
\label{a1}
J[{\bf u}] = h\left(t_\text{f}\hspace{1pt},{\bf x}(t_\text{f}\hspace{1pt})\right) 
+ \int_{0}^{t_\text{f}\hspace{1pt}} {g\left( {\bf x}(t),{\bf u}(t),t \right)} \text{\rm d}t
\end{equation}
subject to the fractional dynamic system
\begin{equation} 
\label{a2}
{}_0^{\text{\tiny C}\hspace{1pt}}{\mathcal D}_t^\alpha  {\bf x}(t) 
= {\bf f}({\bf x}(t),{\bf u}(t),t), \quad 0<\alpha\leq 1,
\end{equation}
with the initial and boundary conditions
\begin{eqnarray}
&& {\bf x}(0) = {\bf x}_0,\label{a3} \\
&& \boldsymbol\psi({\bf x}({t_\text{f}\hspace{1pt}}),t_\text{f}\hspace{1pt}) = \mathbf 0, \label{w4}
\end{eqnarray}
and with the path constraint
\begin{equation} 
\label{w1}
\boldsymbol\phi\left({\bf x}(t),{\bf u}(t),t\right) \le \mathbf 0.
\end{equation}
\end{subequations}
Here, functions $h$, $g$, $\mathbf f$, $\boldsymbol{\psi}$ 
and $\boldsymbol{\phi}$ are sufficiently continuously differentiable 
and defined by the following mappings:
\begin{eqnarray*}
&& h: \mathbb{R}\times \mathbb{R}^p\rightarrow \mathbb{R},\\
&& g: \mathbb{R}^p\times \mathbb{R}^q\times \mathbb{R}\rightarrow \mathbb{R},\\
&& \mathbf f: \mathbb{R}^p\times \mathbb{R}^q\times \mathbb{R}\rightarrow \mathbb{R}^p,\\
&& \boldsymbol{\psi}: \mathbb{R}^p\times \mathbb{R}\rightarrow \mathbb{R}^{r_1}, \ \ \ 0\le r_1\le q,\\
&& \boldsymbol{\phi}: \mathbb{R}^p\times \mathbb{R}^q\times \mathbb{R}\rightarrow
\mathbb{R}^{r_2},\ \ \ 0 \le r_2.
\end{eqnarray*}
Most FOCPs found in the literature, in various disciplines like engineering, control, signal processing 
and others, belong to the above general class of FOCPs 
\cite{ding2012optimal,sweilam2016optimal,sweilam2017legendre}.
 
For easy application of our numerical method, the time domain 
$[0, t_\text{f}\hspace{1pt}]$ is mapped to the canonical interval $[0,1]$, 
by using the affine transformation
$t\rightarrow \tau t_\text{f}\hspace{1pt}$.
Applying this mapping and by noting that
$$
{}_0^C\mathop {\mathcal D}\nolimits_t^\alpha  {\bf x}(t)
=(t_\text{f}\hspace{1pt})^{-\alpha}\  
{}_0^C\mathop {\mathcal D}\nolimits_\tau^\alpha{\bf x}(\tau),
$$
the optimal control problem \eqref{main-focp} is converted 
to the following form:
\begin{subequations}
\label{main-focp2}
\begin{align}[left = \empheqlbrace\,]
\min J[{\bf u}] = &h\left(t_\text{f}\hspace{1pt},{\bf x}(1)\right) 
+ t_\text{f}\int_{0}^{1} {g\left( {\bf x}(\tau),{\bf 
u}(\tau),t_\text{f}\hspace{1pt}\tau
\right)} \text{\rm d}\tau\label{focpJ}\\
\text{s.t.}\ \ &{}_0^{\text{\tiny C}\hspace{1pt}}{\mathcal D}_t^\alpha  {\bf x}(\tau) 
= (t_\text{f}\hspace{1pt})^\alpha {\bf f}({\bf x}(\tau),{\bf u}(\tau),t_\text{f}\hspace{1pt}\tau),
\label{focpD}\\
& {\bf x}(0) = {\bf x}_0, \label{focpIni}\\
& \boldsymbol\psi({\bf x}({1}),t_\text{f}\hspace{1pt}) =\mathbf  0,\label{focpTer}\\
&\boldsymbol\phi({\bf x}(\tau),{\bf u}(\tau),t_\text{f}\hspace{1pt}\tau) \le \mathbf 0.\label{focpPath}
\end{align}
\end{subequations}
It should to be noted that, after applying the mentioned transformation,
the symbols of variables in \eqref{main-focp2} should be changed to new symbols. 
However, for the sake of simplicity, we shall retain the symbols already used.


\section{The proposed method}
\label{sec:4}

We present a direct numerical approach to solve the FOCPs stated in \eqref{main-focp2}.  
Let the interval $[0,1]$ be divided into $n$ equal parts by the following mesh points:
\[
\tau_{k}:=kh,\ k=0,\dots,n,
\]
where $h=1/n$. Moreover, for the sake of simplicity, 
we consider the notations
\begin{align}
\label{notations}
\mathbf x_i:=\mathbf x(\tau_i),
\ \ \ \mathbf u_i:=\mathbf u(\tau_i),\ \ \
\mathbf f_i:=\mathbf f(\mathbf x(\tau_i),
\mathbf u(\tau_i),t_\text{f}\hspace{1pt}\tau_i).
\end{align}


\subsection{Discretization of the performance index}

The performance index \eqref{focpJ} is discretized by a quadrature formula, 
such as the trapezoidal, Simpson or other rules. In this way, we have
\begin{equation} 
\label{p1}
J_{n}= h(t_\text{f}\hspace{1pt}, \mathbf x_n)+ t_\text{f}\hspace{1pt}
\sum_{i=0}^{n}w_{i}g(\mathbf x_{i},\mathbf  u_{i}, t_\text{f}\hspace{1pt}\tau_{i}),
\end{equation}
where, in the case of using the trapezoidal rule, 
\[
w_0:=\frac h 2,\quad w_i:=h,\ i=1,\dots,n-1,\quad w_n:=\frac h 2
\] 
and, in the case of using Simpson rule,
\[
w_0:=\frac h 3,\quad w_{2i-1}:=\frac {4h} 3,\ 
w_{2i}:=\frac {2h} 3,\ i=1,\dots,\frac n 2-1,\quad w_n:=\frac h 3.
\] 


\subsection{Discretization of the fractional dynamic equation}

By applying the fractional integral operator ${}_0\mathcal I^\alpha_\tau$ 
to both sides of the fractional dynamic  equation \eqref{focpD}, 
and by noting that ${}_0\mathcal I^\alpha_\tau \left( {}_0^{\text{\tiny C}
\hspace{1pt}}\mathcal D^\alpha_\tau \mathbf x\left(\tau \right)\right) 
= \mathbf x(\tau) - \mathbf x(0)$,
the above fractional dynamic equations can be converted 
into the following fractional integral form:
\begin{equation}
\mathbf x(\tau)= \mathbf  x(0)+ (t_\text{f}\hspace{1pt})^\alpha\,{}_0\mathcal I_{\tau}^\alpha
\mathbf f(\mathbf x(\tau),\mathbf u(\tau), t_\text{f}\hspace{1pt}\tau).
\end{equation}
By collocating the above equation at $\tau= \tau_{i}$, $i=1, \ldots,n$, we have
\begin{equation}
\label{collocatedeq}
\mathbf x(\tau_{i})= \mathbf x(0) + \left(t_\text{f}\hspace{1pt}\right)^\alpha 
\left[{}_0\mathcal I_{t}^\alpha\left(\mathbf f(\mathbf x(\tau), \mathbf u(\tau),
t_\text{f}\hspace{1pt}\tau)\right)\right]_{\tau=\tau_i}.
\end{equation}
Considering the initial condition \eqref{focpIni} and notations \eqref{notations}, 
and by applying GL/TR/SI approximation formulas 
\eqref{GLformula}/\eqref{TRformula}/\eqref{SIformula}, 
we can approximate the equation \eqref{collocatedeq} as follows:
\begin{equation}
\mathbf x_{i}= \mathbf x_{0} + (t_\text{f}\hspace{1pt})^\alpha\sum_{j=0}^{i} 
\varpi^{(\alpha)}_{ij}\mathbf f(\mathbf x_{j},
\mathbf u_{j}, t_\text{f}\hspace{1pt}\tau_{j}), 
\quad i=1,2,\ldots,n,
\end{equation}
where $\varpi^{(\alpha)}_{ij}$ are the coefficients of GL/TR/SI approximation formula 
defined in \eqref{WoGL}/\eqref{WoTR}/\eqref{WoSI}.
Similarly, we can discretize the path constraint \eqref{focpPath} as follows:
\begin{equation}
\boldsymbol{\phi}(\mathbf x_i,\mathbf u_i,
t_\text{f}\hspace{1pt}\tau_i)\le \mathbf 0,\ \ i=0,\dots,n.
\end{equation}
Moreover, the terminal condition \eqref{focpTer} can be discretized as
\begin{equation}
\boldsymbol{\psi}(\mathbf x_n,t_\text{f}\hspace{1pt})=\mathbf 0.
\end{equation}
In summary, the FOCP \eqref{main-focp2} is transcribed 
into the following nonlinear programming problem (NLP):
\begin{subequations}
\label{main-op}
\begin{align}[left = \empheqlbrace\,]
\min J_n &=  h(t_\text{f}\hspace{1pt}, \mathbf x_n)
+ t_\text{f}\hspace{1pt}\sum_{i=0}^{n}w_{i}g(\mathbf x_{i},\mathbf  u_{i}, 
t_\text{f}\hspace{1pt}\tau_{i}) \label{opJ}\\
\text{s.t.}\ \ & \mathbf x_{i}= \mathbf x_{0} + (t_\text{f}\hspace{1pt})^\alpha
\sum_{j=0}^{i} \varpi^{(\alpha)}_{ij}\mathbf f\left(\mathbf x_{j},
\mathbf u_{j}, t_\text{f}\hspace{1pt}	\tau_{j}\right), 
\quad i=0,1,\ldots ,n	,\label{opD}\\
& \boldsymbol\psi(\mathbf x_n,t_\text{f}\hspace{1pt}) =\mathbf  0,\label{oTer}\\
&\boldsymbol\phi(\mathbf x_i,\mathbf  u_i,t_\text{f}\hspace{1pt}\tau_i) 
\le \mathbf 0,\ \ i=0,\dots,n.\label{oPath}
\end{align}
\end{subequations}
It should be remembered that the decision variables of the above optimization problem are:
$\mathbf x_i$, $i = 0,\dots,n$, $\mathbf u_i$, $i = 0,\dots,n$, and maybe $t_\text{f}\hspace{1pt}$. 
By utilizing an optimization solver, we can solve the optimization problem 
\eqref{main-op} and find the optimal value of the decision variables, which 
are the approximations of state and control functions 
at the points $\tau_i$, $i = 0,\dots, n$.


\subsubsection{Reformat of the resulted optimization problem to the classical form}

Traditionally, in an NLP, it is preferred that the decision variables 
are placed in a vector. This helps to analyze and utilize a solver for solving the problem.  
Here, we reformulate the NLP \eqref{main-op} into the standard form, where the decision 
variables are collected in a vector $\mathbf z$. For fixed final time problems, 
the vector of all decision variables is
\begin{equation*}
\mathbf z:=\left[\mathbf x_0,\mathbf u_0,\mathbf x_1,
\mathbf u_1,\dots,\mathbf x_n,\mathbf u_n \right]^T
\end{equation*}
and for free final time problems
\begin{equation*}
\mathbf z:=\left[\mathbf x_0,\mathbf u_0,\mathbf x_1,\mathbf u_1,\dots,
\mathbf x_n,\mathbf u_n,t_\text{f}\hspace{1pt} \right]^T.
\end{equation*}
If we define new functions $\widehat h$ and $\widehat g$  as
\begin{equation}
\widehat h(\mathbf z):=h(t_\text{f}\hspace{1pt}, \mathbf x_n),\ \ \
\mathbf {\widehat g}(\mathbf z):=
\begin{bmatrix}
g(\mathbf x_0, \mathbf u_0, t_\text{f}\hspace{1pt}\tau_0)\\
g(\mathbf x_1, \mathbf u_1, t_\text{f}\hspace{1pt}\tau_1)\\
\vdots\\
g(\mathbf x_n, \mathbf u_n, t_\text{f}\hspace{1pt}\tau_n)
\end{bmatrix},
\end{equation}
then we can express the objective function \eqref{opJ} as
\[
J_n (\mathbf z):= \widehat h(\mathbf z)+ t_\text{f}\hspace{1pt} 
\mathbf{w}^{T}\mathbf{\widehat g}(\mathbf z),
\]
where $\mathbf w$ is the vector of quadrature weights in \eqref{p1}, i.e.,
\begin{equation}
\mathbf w:=\left[w_0,\dots,w_n\right]^T.
\end{equation}
It should be noted that, in case of using trapezoidal quadrature, we have 
\begin{equation}
\label{wTR}
\mathbf w=\mathbf w_{\hspace{-1pt}\text{\tiny TR}}
:=\left[\tfrac h 2,h,\dots,h,\tfrac h 2\right]
\end{equation}
and, for Simpson quadrature, we have
\begin{equation}
\label{wSI}
\mathbf w=\mathbf w_{\hspace{-1pt}\text{\tiny SI}}:=\left[\tfrac h 3,\tfrac {4h} 3,
\tfrac {2h}3,\tfrac {4h} 3,\tfrac {2h}3,\dots,\tfrac {4h} 3,\tfrac {2h}3,\tfrac h 3\right].
\end{equation}
The constraints \eqref{opD}, for $i=0,1,\dots,n$, can be represented 
in the following matrix form:
\begin{equation}
\label{Dmatf}
\mathbf {X}= \mathbf {X}_{0} + (t_\text{f}\hspace{1pt})^\alpha\hspace{.2mm} 
\mathbf  F\hspace{1pt} \left[\mathbf  W^{(\alpha)}\right]^T,
\end{equation}
where $\mathbf X$, $\mathbf X_0$ and $\mathbf F$ are 
$p\times(n+1)$ matrices defined as
\begin{equation}
\label{ggg}
\mathbf X  =
\begin{bmatrix}
\mathbf x_0&\dots&\mathbf x_n
\end{bmatrix},
\quad
\mathbf X_0=\begin{bmatrix}
\mathbf x_0&\dots&\mathbf x_0
\end{bmatrix},
\quad  
\mathbf F =
\begin{bmatrix}
\mathbf f(\mathbf  x_0,\mathbf u_0, t_\text{f}\hspace{1pt}\tau_0)
&\dots& \mathbf f(\mathbf  x_n,\mathbf u_n,t_\text{f}\hspace{1pt}\tau_n)\\
\end{bmatrix},
\end{equation}
and $\mathbf W^{(\alpha)}$ is the fractional GL/TR/SI integration matrix 
defined in \eqref{GLmat}/\eqref{TRmat}/\eqref{SImat}.
The  matrix equation \eqref{Dmatf} can be reformulated as follows:
\begin{equation}
\textbf{vec}\left(\mathbf {X}\right)= \textbf{vec}\left(\mathbf {X}_{0}\right) 
+ t_\text{f}\hspace{1pt}^\alpha \textbf{vec}\left( \mathbf  
F \left[\mathbf  W^{(\alpha)}\right]^T  \right),\label{MatEq}
\end{equation}
where $\textbf{vec}$ is the vectorization operator, which converts 
a matrix into a column vector by stacking the columns 
of the matrix on the top of each other.
  
Let $\mathbf I_{k}$ be the identity matrix of order $k$ 
and the matrix $\mathbf R$ be defined as
\begin{equation*}
\mathbf R:=
\begin{cases}
\left[\mathbf I_p \ |\ \mathbf 0_{p\times q}\right],
& \text{if $t_\text{f}\hspace{1pt}^\alpha$ be fixed},\\
\left[\mathbf I_p \ |\ \mathbf 0_{p\times q}\ |\ 0\right],
& \text{if $t_\text{f}\hspace{1pt}^\alpha$ be free}.	
\end{cases}
\end{equation*} 
Then, we can express $\textbf{vec}\left( \mathbf X\right)$  
based on $\mathbf z$ as
\begin{equation}
\label{vecx}
\textbf{vec}\left( \mathbf X\right)
=\left(\mathbf I_n\otimes \mathbf R \right) \mathbf z,
\end{equation} 
where $\otimes$ denotes the Kronecker product (or tensor product). 
Moreover, if $\mathbf 1_n$ is a row $n$-vector whose entries are all 
equal to $1$, then we can express $\text{vec}(\mathbf X_0)$ based on $\mathbf z$ as 
\begin{equation}
\label{vecx0}
\textbf{vec}(\mathbf X_0)=\mathbf x_0\otimes \mathbf 1_n.
\end{equation}
To express $\textbf{vec}\left(  \mathbf  F \left[\mathbf  W^{(\alpha)}\right]^T\right)$ 
based on $\mathbf z$, we use Theorem~13.26 in \cite{laub2005matrix}, which states 
that for any three matrices $\mathbf A$, $\mathbf B$, and $\mathbf E$ 
for which the matrix product $\mathbf {ABE}$ is defined, one has 
\begin{equation*}
\textbf{vec}(\mathbf {ABE}) 
= (\mathbf E^T\otimes \mathbf A)\textbf{vec}(\mathbf B).
\end{equation*}
Now, by using the above equation and by noting that 
$\mathbf  F \left[\mathbf  W^{(\alpha)}\right]^T
=\mathbf I_p \mathbf  F   \left[\mathbf  W^{(\alpha)}\right]^T$,
 we conclude that
\begin{equation}
\label{vecfw}
\textbf{vec}\left( \mathbf  F  \left[\mathbf  W^{(\alpha)}\right]^T\right)
=\left(\mathbf W^{(\alpha)} \otimes \mathbf I_p\right) \textbf{vec}\left( \mathbf  F\right).
\end{equation}  
In view of \eqref{ggg}, we can write 
\begin{equation}
\label{vecfz}
\textbf{vec}\left( \mathbf  F\right)=\mathbf {\widehat f}(\mathbf z):=
\begin{bmatrix}
\mathbf f(\mathbf x_0, \mathbf u_0, t_\text{f}\hspace{1pt}\tau_0)\\
\mathbf f(\mathbf x_1, \mathbf u_1, t_\text{f}\hspace{1pt}\tau_1)\\
\vdots\\
\mathbf f(\mathbf x_n, \mathbf u_n, t_\text{f}\hspace{1pt}\tau_n)
\end{bmatrix}.
\end{equation}
Using the above notation and equations \eqref{vecx}--\eqref{vecfw},   
we can write \eqref{MatEq} as
\[
\mathbf c(\mathbf z):= \left(\mathbf I_n\otimes \mathbf R \right) 
\mathbf {z}- \mathbf x_0\otimes \mathbf 1_n 
-(t_\text{f}\hspace{1pt})^\alpha \left(\mathbf W^{(\alpha)} \otimes \mathbf I_p\right) 
\mathbf {\widehat f}(\mathbf z)=\mathbf 0.
\]
To replace the constraints \eqref{oTer} and \eqref{oPath} by other constraints 
based on $\mathbf z$, we define functions $\boldsymbol {\widehat \phi} $ 
and $\boldsymbol {\widehat \psi}$ as
\begin{equation*}
\boldsymbol {\widehat \phi}(\mathbf z)=
\begin{bmatrix}
\boldsymbol { \phi}(\mathbf x_0, \mathbf u_0, t_\text{f}\hspace{1pt}\tau_0)\\
\boldsymbol { \phi}(\mathbf x_1, \mathbf u_1, t_\text{f}\hspace{1pt}\tau_1)\\
\vdots\\
\boldsymbol { \phi}(\mathbf x_n, \mathbf u_n, t_\text{f}\hspace{1pt}\tau_n)
\end{bmatrix},
\qquad
\boldsymbol {\widehat \psi}(\mathbf z)=
\boldsymbol { \psi}(\mathbf x_n,t_\text{f}\hspace{1pt}).
\end{equation*}     
Then the constraints \eqref{oTer} and \eqref{oPath} can be expressed as
\begin{equation*}
\boldsymbol{\widehat \psi}(\mathbf z) =\mathbf  0,
\qquad \boldsymbol{\widehat\phi}(\mathbf z) \le \mathbf 0.
\end{equation*}
The above results are summarized in the following theorem.

\begin{thm}
The optimization problem \eqref{main-op} can be expressed 
in the standard form
\begin{subequations}
\label{main-Sop}
\begin{align}[left = \empheqlbrace\,]
\min \widehat J_n (\mathbf z)= &\widehat h(\mathbf z)
+ t_\text{\rm f}\hspace{1pt} \mathbf{w}^{T}\mathbf{\widehat g}(\mathbf z) \label{SopJ}\\
\text{s.t.}\ \ &\mathbf c(\mathbf z):=\left(\mathbf I_n\otimes \mathbf R \right) \mathbf{z}
- \mathbf x_0\otimes \mathbf 1_n - (t_\text{\rm f}\hspace{1pt})^\alpha \left(\mathbf W^{(\alpha)} 
\otimes \mathbf I_p\right) \mathbf {\widehat f}(\mathbf z)=\mathbf 0,\label{SopD}\\
& \boldsymbol{\widehat \psi}(\mathbf z) =\mathbf  0,\label{SoTer}\\
&\boldsymbol{\widehat\phi}(\mathbf z) \le \mathbf 0.\label{SoPath}
\end{align}
\end{subequations}	
\end{thm}

In summary, the solution of the FOCP \eqref{main-focp} is reduced to the solution 
of the NLP \eqref{main-Sop}. In this NLP, $\mathbf w$ can be 
$\mathbf w_{\hspace{-1pt}\text{\tiny TR}}$ or $\mathbf w_{\hspace{-1pt}\text{\tiny SI}}$ 
defined in \eqref{wTR} and \eqref{wSI}, respectively. Moreover, the fractional integration matrix 
$\mathbf W^{(\alpha)}$ can be $\mathbf W^{(\alpha)}_\text{\tiny GL}$, 
$\mathbf W^{(\alpha)}_\text{\tiny TR}$ or $\mathbf W^{(\alpha)}_\text{\tiny SI}$, 
defined in \eqref{GLmat}, \eqref{TRmat} and \eqref{SImat}, respectively.
We categorize our methods as follows:
\begin{itemize}
\item If $\mathbf w=\mathbf w_{\hspace{-1pt}\text{\tiny GL}}$ 
and $\mathbf W^{(\alpha)}=\mathbf W^{(\alpha)}_\text{\tiny GL}$, 
then the method is called the ``direct GL method''.

\item If $\mathbf w=\mathbf w_{\hspace{-1pt}\text{\tiny TR}}$ and 
$\mathbf W^{(\alpha)}=\mathbf W^{(\alpha)}_\text{\tiny TR}$, 
then the method is called the ``direct TR method''.

\item If $\mathbf w=\mathbf w_{\hspace{-1pt}\text{\tiny SI}}$ and 
$\mathbf W^{(\alpha)}=\mathbf W^{(\alpha)}_\text{\tiny SI}$, 
then the method is called the ``direct SI method''.
\end{itemize}  


\subsubsection{Gradient of the objective function and Jacobian of the constraints}

The speed of the proposed method depends on the speed one solves the NLP \eqref{main-Sop}. 
On the other hand, a crucial task in solving NLPs  is computing the first derivatives 
of the objective and constraint functions (gradients of the objective function 
and Jacobians of the constraints).
Many NLP-solvers allow users to supply the exact first derivatives. 
If the derivatives are not provided, then the solvers approximate 
them numerically, by finite difference formulas.
However, the use of first derivative approximations may seriously degrade 
the performance and convergence of the NLP-solver. Thus, it is highly recommended 
that at least the exact gradient of the objective function and the Jacobian 
of the constraints are provided by the user. These can make a great 
improvement in the convergence, accuracy and computation time of the NLP-solver.
In what follows, we detail the closed forms for the 
gradient of the objective function and the Jacobian of constraints.

For free time problems,	the gradient of the objective function 
can be obtained as
\begin{equation}
\nabla \widehat J_n (\mathbf z)=
\begin{bmatrix}
\partial \widehat J_n (\mathbf z) / \partial \mathbf x_0\\
\partial \widehat J_n (\mathbf z) / \partial \mathbf u_0\\
\partial \widehat J_n (\mathbf z) / \partial \mathbf x_1\\
\partial \widehat J_n (\mathbf z) / \partial \mathbf u_1\\
\vdots\\
\partial \widehat J_n (\mathbf z) / \partial \mathbf x_n\\
\partial \widehat J_n (\mathbf z) / \partial \mathbf u_n\\
\partial \widehat J_n (\mathbf z) / \partial  t_\text{f}\hspace{1pt}\\
\end{bmatrix}
=
\begin{bmatrix}
t_\text{f}\hspace{1pt}\, w_0\, g_{\mathbf x}(\mathbf x_0,\mathbf u_0,t_\text{f}\hspace{1pt}\tau_0)\\
t_\text{f}\hspace{1pt}\, w_0\, g_{\mathbf u}(\mathbf x_0,\mathbf u_0,t_\text{f}\hspace{1pt}\tau_0)\\
t_\text{f}\hspace{1pt}\, w_1\, g_{\mathbf x}(\mathbf x_1,\mathbf u_1,t_\text{f}\hspace{1pt}\tau_1)\\
t_\text{f}\hspace{1pt}\, w_1\, g_{\mathbf u}(\mathbf x_1,\mathbf u_1,t_\text{f}\hspace{1pt}\tau_1)\\
\vdots\\
h_{\mathbf x}(t_\text{f}\hspace{1pt},\mathbf x_n)+t_\text{f}\hspace{1pt} w_n g_{\mathbf x}(\mathbf x_n,
\mathbf u_n,t_\text{f}\hspace{1pt}\tau_n)\\
t_\text{f}\hspace{1pt} w_n g_{\mathbf u}(\mathbf x_n,\mathbf u_n,t_\text{f}\hspace{1pt}\tau_n)\\
h_{t}(t_\text{f}\hspace{1pt},\mathbf x_n)+ \sum_{i=0}^n w_i \big[g(\mathbf x_i,\mathbf
u_i,t_\text{f}\hspace{1pt}\tau_i)+t_\text{f}\hspace{1pt} \tau_i\,
g_t(\mathbf x_i,\mathbf u_i,t_\text{f}\hspace{1pt}\tau_i)\big]
\end{bmatrix}.
\end{equation}
For fixed final time problems, the gradient is the same 
as above, except that the last row is removed.

The Jacobian of constraints \eqref{SopD}, 
for free final time problems, can be derived as
\begin{equation}
\nabla \mathbf c(\mathbf z)=\left[
\begin{array}{c|c}
\left(\mathbf I_n\otimes \mathbf R \right) - \left(t_\text{f}\hspace{1pt}\right)^\alpha 
\left(\mathbf W^T \otimes \mathbf I_p\right)\nabla\mathbf{\widehat f}\left(\mathbf z\right) 
& (t_\text{f}\hspace{1pt})^{\alpha-1} \left(\mathbf W^T \otimes \mathbf I_p\right)\left(
-\alpha\mathbf{\widehat f}\left(\mathbf z\right)
- t_\text{f}\hspace{1pt}\mathbf{\widehat f}_{t_\text{f}\hspace{1pt}}\left(\mathbf z\right)\right)
\end{array}
\right],
\end{equation}
such that
\[
\nabla\mathbf{\widehat 	f}\left(\mathbf z\right):=
\begin{bmatrix}
[\mathbf f_{\mathbf x}]_0 & [\mathbf f_{\mathbf u}]_0 && &&&\\
&& [\mathbf f_{\mathbf x}]_1 & [\mathbf f_{\mathbf u}]_1&& &\\
&&& &\ddots&&\\
&&&& & [\mathbf f_{\mathbf x}]_n & [\mathbf f_{\mathbf u}]_n
\end{bmatrix},
\quad
\mathbf{\widehat f}_{t_\text{f}\hspace{1pt}}\left(\mathbf z\right):=
\begin{bmatrix}
\tau_0\mathbf f_t(\mathbf x_0, \mathbf u_0, t_\text{f}\hspace{1pt}\tau_0)\\
\tau_1\mathbf f_t(\mathbf x_1, \mathbf u_1, t_\text{f}\hspace{1pt}\tau_1)\\
\vdots\\
\tau_n\mathbf f_t(\mathbf x_n, \mathbf u_n, t_\text{f}\hspace{1pt}\tau_n)
\end{bmatrix},
\]
where $[\mathbf f_{\mathbf x}]_i$ and $[\mathbf f_{\mathbf u}]_i$ 
are the $p$-vectors defined as
\[
[\mathbf f_{\mathbf x}]_i:=\mathbf f_{\mathbf x}(\mathbf x_i,
\mathbf u_i,t_\text{f}\hspace{1pt}\tau_i), \ \
[\mathbf f_{\mathbf u}]_i:=\mathbf f_{\mathbf u}(\mathbf x_i,
\mathbf u_i,t_\text{f}\hspace{1pt}\tau_i).
\]
In the case of fixed final time problems, 
the Jacobian of constraints \eqref{SopD} is obtained as
\begin{equation}
\nabla \mathbf c(\mathbf z)=\left[
\begin{array}{c}
\left(\mathbf I_n\otimes \mathbf R \right)  
- (t_\text{f}\hspace{1pt})^\alpha \left(\mathbf W^T 
\otimes \mathbf I_p\right)\nabla\mathbf{\widehat f}\left(\mathbf z\right) 
\end{array}
\right].
\end{equation}
The Jacobian of constraint \eqref{SoTer}, for fixed final time FOCPs, 
is derived as
\begin{equation}
\nabla  \boldsymbol{\widehat \psi}(\mathbf z) =\left[
\begin{array}{c|c|c}
\mathbf 0_{r_1\times(n-1)(p+q)} 
& \boldsymbol{ \psi}_{\mathbf x}(\mathbf x_n,t_\text{f}\hspace{1pt}) 
& \mathbf 0_{r_1\times q}
\end{array}
\right]
\end{equation} 
while for free final time FOCPs is given by
\begin{equation}
\nabla  \boldsymbol{\widehat \psi}(\mathbf z) =\left[
\begin{array}{c|c|c|c}
\mathbf 0_{r_1\times(n-1)(p+q)} 
& \boldsymbol{ \psi}_{\mathbf x}(\mathbf x_n,t_\text{f}\hspace{1pt}) 
& \mathbf 0_{r_1\times q}
&\boldsymbol{ \psi}_{t}(\mathbf x_n,t_\text{f}\hspace{1pt})
\end{array}
\right].
\end{equation}
In case of free final time FOCPs, the Jacobian of the inequality 
constraints \eqref{SoPath} is obtained as
\begin{equation}
\label{IaT}
\nabla \boldsymbol{\widehat \phi}(\mathbf z) 
=
\begin{bmatrix}
[\boldsymbol{\phi}_{\mathbf x}]_0
& [\boldsymbol{\phi}_{\mathbf u}]_0 
&& &&&&\tau_0 \boldsymbol{\phi}_{t}(\mathbf x_0,\mathbf u_0,t_\text{f}\hspace{1pt}\tau_0)\\
&& [\boldsymbol{\phi}_{\mathbf x}]_1 & [\boldsymbol{\phi}_{\mathbf u}]_1 
&& &&\tau_1 \boldsymbol{\phi}_{t}(\mathbf x_1,\mathbf u_1,t_\text{f}\hspace{1pt}\tau_1)\\
&&& &\ddots&&\\
&&&& & [\boldsymbol{\phi}_{\mathbf x}]_n &[\boldsymbol{\phi}_{\mathbf u}]_n
&\tau_n \boldsymbol{\phi}_{t}(\mathbf x_n,\mathbf u_n,t_\text{f}\hspace{1pt}\tau_n)
\end{bmatrix},
\end{equation}
where $[\boldsymbol{\phi}_{\mathbf x}]_i$ and $[\boldsymbol{\phi}_{\mathbf u}]_i$ 
are the $r_2$-vectors defined by
\[
[\boldsymbol{\phi}_{\mathbf x}]_i:=\boldsymbol{\phi}_{\mathbf x}(\mathbf x_i,
\mathbf u_i,t_\text{f}\hspace{1pt}\tau_i),\quad
[\boldsymbol{\phi}_{\mathbf u}]_i:=\boldsymbol{\phi}_{\mathbf u}(\mathbf x_i,
\mathbf u_i,t_\text{f}\hspace{1pt}\tau_i), \quad i=0,\dots,n.
\]
In case of fixed final time problems, the last column of the 
Jacobian matrix \eqref{IaT} is removed.

As we can see, the partial derivatives of $g$, $\mathbf f$, 
$\boldsymbol{ \psi}$ and $\boldsymbol{ \phi}$
with respect to $\mathbf x$ and $\mathbf u$ are needed to
provide the gradient of the objective function and the Jacobian of the constraints.
We can use symbolic computation to supply these partial derivatives.


\section{Numerical examples}
\label{sec:5}

We now illustrate the direct methods presented in Section~\ref{sec:4} 
through numerical experiments. We have implemented the direct GL, TR and SI methods 
using \textsc{Matlab} on a 3.5 GHz Core i7 personal computer with 8 GB of RAM.
Moreover, for  solving  the nonlinear programming problem \eqref{main-Sop}, 
the solver \textsc{Ipopt} \cite{Waechter2006} was used, which is based on 
an interior-point algorithm. In \textsc{Ipopt}, we can adjust the accuracy 
of solution by the input parameter \texttt{tolrfun}. In our numerical experiments, 
we set \texttt{tolrfun}=$10^{-12}$.
	
We consider four nontrivial examples. The first two examples are new while 
the last two  examples have been investigated before in the literature. 
At each example, the capability of the proposed direct methods is highlighted.
For the first example, which has an exact solution, the accuracy 
of our methods, and the required CPU time, are assessed.
In the second example, the ability of the method in solving a free final 
time FOCP with path constraints is investigated.
In the third example, we show that the proposed methods are also
suitable to deal with bang-bang optimal control problems. 
In the last example, we assess the ability of our direct methods 
in solving a practical and challenging problem of HIV.


\subsection{Example 1: A FOCP with exact solution}

In this example, we consider the following nonlinear FOCP:
\begin{equation}
\label{eq:ex1}
\begin{aligned}
\min J[{\bf u}] &= \int_{0}^{20} { \left[ 1- \big( x-0.01\,{t}^{2}-1 \big) ^{2}
+u-2\,\sqrt{\pi }{{\sl J}_{0}\big(4\,\sqrt {t}\big)} \right] ^{2}}\, \text{d}t\\
\text{s.t.}\ \ &{}_0^{\text{\tiny C}\hspace{1pt}}{\mathcal D}_t^\alpha  { x}(t) 
= - \big( x(t)-0.01\,{t}^{2}-1 \big) ^{2}+u(t)+1+{\frac {2\,{t}^{3/2}}{75\,\sqrt {\pi }}},\\
& {x}(0) = 1,\\
& x(20)=5+\sin(8\sqrt 5).
\end{aligned}
\end{equation}
The exact solution for $\alpha=1/2$ is 
\begin{align*}
&u_{\text{ex}}(t)=	- \cos^2 \big( 4\,\sqrt {t} \big)  
+2\,\sqrt {\pi}{{\sl J}_{0}\big(4\,\sqrt {t}\big)},\\
&x_{\text{ex}}(t)=\sin \big( 4\,\sqrt {t} \big) +0.01\, {t}^2+1,
\end{align*}
where $J_0$ is the first kind Bessel function of order zero \cite{abramowitz1964handbook}. 
The domain of this problem is large and the control and state functions have 
algebraic singularities and an oscillating behavior. Consequently, this example is a suitable 
test problem for assessing the accuracy and efficiency of the presented methods. 

\begin{figure}
\includegraphics[width=\textwidth]{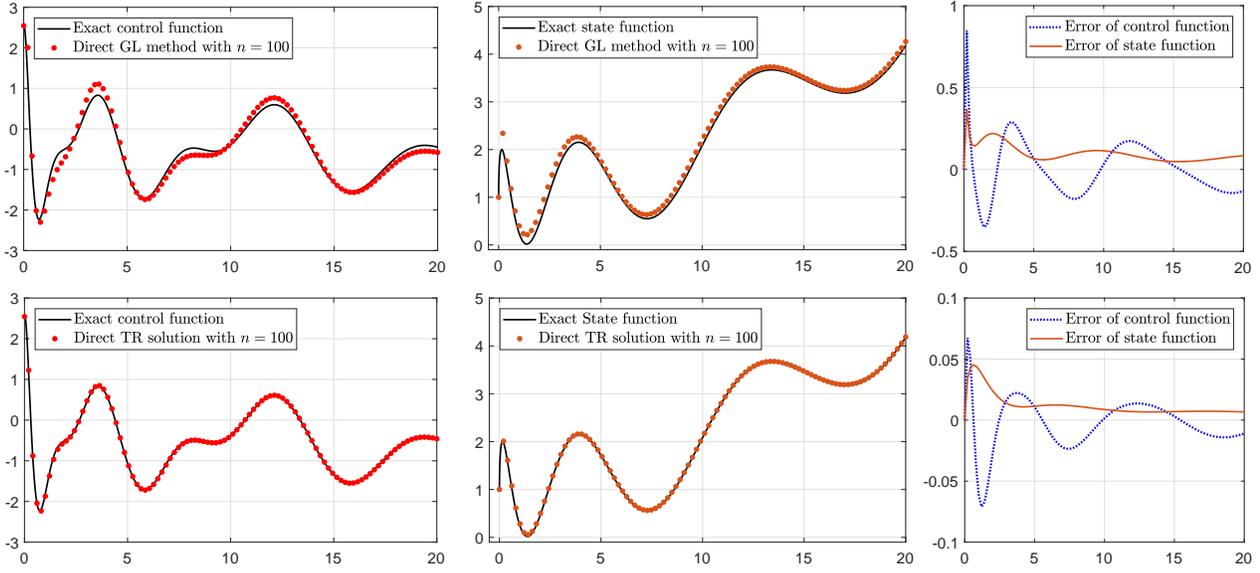}
\caption{(Example 1, with $\alpha=0.5$) Comparison of the exact and obtained 
solutions and  error functions. Above: Direct GL method with $n=100$.  
Below: Direct TR method with $n=100$.}
\label{FIG1}
\end{figure}

\begin{figure}[ht]
\includegraphics[width=\textwidth]{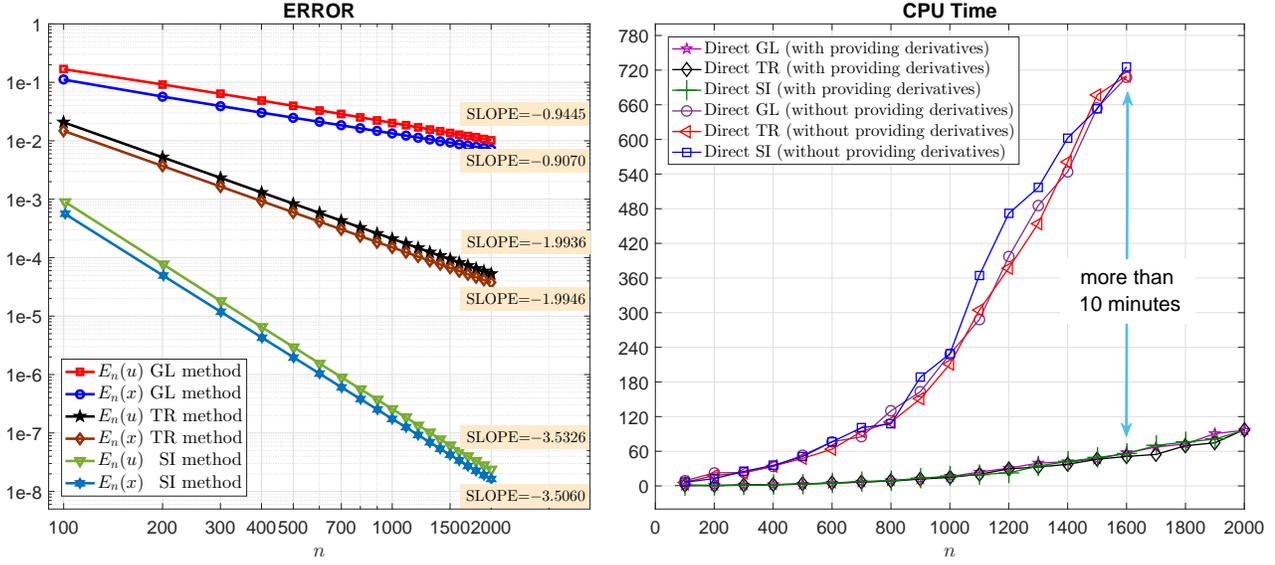}
\caption {(Example 1, with $\alpha=0.5$) Left: The measured 
$E_n(u)$ and $E_n(x)$ in the direct GL/TR/SL methods  
against various $n$ in \texttt{loglog} scale.
Right: CPU time of direct GL/TR/SI methods, with 
and without providing derivatives, versus $n$.}
\label{FIG2}
\end{figure}

By applying the direct GL and TR methods with $n=100$, the obtained control, state and error functions 
are plotted in Figure~\ref{FIG1}. We see that, with $n=100$, the accuracy of the direct GL 
method is not satisfactory. However, the direct TR method generates a solution with reasonable accuracy. 
The results of the direct GL, TR  and SI methods, with various values of $n$, 
are reported in Table~\ref{TAB1}. In this table, to show the impact of using 
the gradient of the objective function and the Jacobian of the constraints, we compare 
the CPU times in the two cases ``NPD'' and ``PD''. By ``NPD'' we refer to the CPU 
time of the methods without providing derivatives of the objective and constraint functions. 
In contrast, ``PD'' refers to the CPU time when the first derivatives of the objective 
and constraint functions are supplied to the direct methods.
In addition, in Table~\ref{TAB1}, $E_n(u)$ and $E_n(x)$ are the $\ell_2$ norm 
of the error in control and state, i.e.,
\[
E_n(u):=\left(\frac 1 n\sum_{i=1}^n(u_i-u_{\rm ex}(\tau_i))^2\right)^{\frac 1 2}, 
\qquad E_n(x):=\left(\frac 1 n\sum_{i=1}^n(x_i-x_{\rm ex}(\tau_i))^2\right)^{\frac 1 2},
\]
where $u_i$ and $x_i$ are the obtained control and state functions in $\tau=\tau_i$ 
by the presented methods with $n$ nodes.
\begin{table}
\caption{(Example~1, with $\alpha=0.5$) Obtained CPU times and the $\ell^2$ norm 
of the errors by direct GL, TR and SI methods for various values of $n$. 
All CPU times are in seconds.}
\label{TAB1}
\resizebox{\textwidth}{!}
{
\begin{tabular}{rlrlrllllllrrrrlllrrrllll}
\hline
&  & \multicolumn{7}{c}{Direct GL method} &  & \multicolumn{7}{c}{Direct TR method} & 
\multicolumn{1}{c}{} & \multicolumn{7}{c}{Direct SI method} \\ 
\cline{3-9}\cline{11-17}\cline{19-25}
\multicolumn{1}{c}{$n$} & \multicolumn{1}{c}{} & \multicolumn{3}{c}{CPU time (s)} & 
\multicolumn{1}{c}{} & \multicolumn{3}{c}{ERROR} & \multicolumn{1}{c}{} & \multicolumn{3}{c}{CPU 
time (s)} & \multicolumn{1}{c}{} & \multicolumn{3}{c}{ERROR} & \multicolumn{1}{c}{} & 
\multicolumn{3}{c}{CPU time (s)} & \multicolumn{1}{c}{} & \multicolumn{3}{c}{ERROR} \\ 
\cline{1-1}\cline{3-5}\cline{7-9}\cline{11-13}\cline{15-17}\cline{19-21}\cline{23-25}
\multicolumn{1}{c}{} & \multicolumn{1}{c}{} & \multicolumn{1}{c}{NPD} & \multicolumn{1}{c}{} & 
\multicolumn{1}{c}{PD} & \multicolumn{1}{c}{} & \multicolumn{1}{c}{$E_n(u)$} & \multicolumn{1}{c}{} 
& \multicolumn{1}{c}{$E_n(x)$} & \multicolumn{1}{c}{} & \multicolumn{1}{c}{NPD} & 
\multicolumn{1}{c}{} & \multicolumn{1}{c}{PD} & \multicolumn{1}{c}{} & \multicolumn{1}{c}{$E_n(u)$} 
& \multicolumn{1}{c}{} & \multicolumn{1}{c}{$E_n(x)$} & \multicolumn{1}{c}{} & 
\multicolumn{1}{c}{NPD} & \multicolumn{1}{c}{} & \multicolumn{1}{c}{PD} & \multicolumn{1}{c}{} 
& \multicolumn{1}{c}{$E_n(u)$} & \multicolumn{1}{c}{} & \multicolumn{1}{c}{$E_n(x)$} \\ 
\cline{1-1}\cline{3-3}\cline{5-5}\cline{7-7}\cline{9-9}\cline{11-11}\cline{13-13}\cline{15-15}
\cline{17-17}\cline{19-19}\cline{21-21}\cline{23-23}\cline{25-25}
100   &&  8.9      && 0.5 && 1.68e-1&&1.11e-1&&6.3       && 0.6 &&2.07e-2 &&1.48e-2
&& 6.4        &&0.5  &&8.99e-4 &&5.60e-4\\
200   && 22.6     && 0.7 && 9.19e-2&&5.71e-2&&17.9     && 0.8 &&5.21e-3 &&3.71e-3
&& 12.7     &&0.8  &&7.66e-5  &&4.91e-5\\
300   && 23.2     && 1.3 && 6.37e-2&&3.94e-2&&22.7     && 1.5 &&2.32e-3 &&1.65e-3
&& 25.2     &&1.1  &&1.80e-5  &&1.18e-5\\
400   && 33.6     && 2.1 && 4.88e-2&&3.04e-2&&34.4     && 2.2 &&1.31e-3 &&9.31e-4
&& 35.6     &&2.1  &&6.48e-6  &&4.30e-6\\
500   && 53.1     && 3.0 && 3.95e-2&&2.48e-2&&47.6     && 3.6 &&8.39e-4 &&5.96e-4
&& 50.8     &&3.2  &&2.94e-6  &&1.97e-6\\
600   && 76.5     && 4.4 && 3.32e-2&&2.11e-2&&63.6     && 5.0 &&5.84e-4 &&4.15e-4
&& 76.5     &&5.1  &&1.54e-6  &&1.04e-6\\
700   && 84.9     &&6.5  && 2.87e-2&&1.84e-2&&94.3     && 6.9 &&4.29e-4 &&3.05e-4
&&100.6    &&7.7  &&8.98e-7  &&6.04e-7\\
800   && 130.1   &&8.7  && 2.52e-2&&1.63e-2&&110.4   && 8.7 &&3.29e-4 &&2.34e-4
&&107.5    &&9.2  &&5.62e-7  &&3.80e-7\\
900   && 163.2   &&12.3&& 2.25e-2&&1.47e-2&&150.9   &&11.9&&2.60e-4 &&1.85e-4
&&188.2    &&13.4&&3.70e-7 &&2.50e-7\\
1000 && 228.1   &&16.1&& 2.03e-2&&1.34e-2&&210.8   &&14.7&&2.11e-4 &&1.50e-4
&&228.9    &&17.0&&2.56e-7 &&1.73e-7\\
1100 && 287.9   &&23.9&& 1.85e-2&&1.23e-2&&304.3   &&20.0&&1.74e-4 &&1.24e-4
&&364.6    &&19.7&&1.83e-7 &&1.24e-7\\
1200 && 397.0   &&30.8&& 1.70e-2&&1.14e-2&&377.3   &&29.6&&1.46e-4 &&1.04e-4
&&471.7    &&22.6&&1.35e-7 &&9.15e-8\\
1300 && 485.5   &&39.4&& 1.57e-2&&1.06e-2&&454.4   &&32.6&&1.25e-4 &&8.87e-5
&&516.7    &&34.6&&1.02e-7 &&6.91e-8\\
1400 && 543.9   &&42.0&& 1.46e-2&&9.89e-3&&561.0   &&37.4&&1.08e-4 &&7.65e-5
&&602.4    &&42.3&&7.84e-8 &&5.34e-8\\
1500 && 654.4   &&46.7&& 1.36e-2&&9.29e-3&&676.9   &&46.2&&9.38e-5 &&6.67e-5
&&652.8    &&49.1&&6.15e-8 &&4.20e-8\\
1600 && 707.8   &&57.3&& 1.28e-2&&8.77e-3&& 710.2  &&51.0&&8.30e-5 &&5.90e-5
&&725.6    &&55.6&&4.92e-8 &&3.36e-7\\
1700 &&$>$800 &&67.2&& 1.20e-2&&8.31e-3&&$>$800&&54.5&&7.30e-5 &&5.19e-5&&$>$800
&&70.2&&4.00e-8 &&2.73e-8\\
1800 &&$>$800 &&72.7&& 1.14e-2&&7.89e-3&&$>$800&&69.7&&6.52e-5 &&4.63e-5&&$>$800
&&76.0&&3.31e-8 &&2.26e-8\\
1900 &&$>$800 &&90.5&& 1.08e-2&&7.52e-3&&$>$800&&74.5&&5.84e-5 &&4.15e-5&&$>$800
&&81.3&&2.79e-8 &&1.90e-6\\
2000 &&$>$800 &&96.4&& 1.03e-2&&7.18e-3&&$>$800&&97.6&&5.26e-5 &&3.74e-5&&$>$800
&&96.6&&2.37e-8 &&1.61e-8\\ \hline
\end{tabular}
}	
\end{table}
Figure~\ref{FIG2} visualizes the results in Table~\ref{TAB1}.  
In Figure~\ref{FIG2} (left), the \texttt{loglog} plot of $E_n(u)$ 
and $E_n(x)$ versus $n$ are plotted. Moreover, the regression lines 
and their slopes are reported too.  In Figure~\ref{FIG2} (right), 
the CPU times of the direct GL, TR and SI methods are plotted 
for the two cases NPD and PD. 
From Table~\ref{TAB1} and Figure~\ref{FIG2}, we can see that 
by supplying the gradient of the objective function and the Jacobian 
of the constraints, the CPU time is significantly reduced. Moreover, 
we find out that the computational times for direct GL, TR and SI methods 
are almost the same. However, the  direct SI method is more accurate 
than the direct TR and GL methods. In view of Figure~\ref{FIG2} (left), 
the accuracy order of the direct GL, TR and SI methods, in solving the fractional
optimal control problem \eqref{eq:ex1}, is $\mathcal O(h)$, $\mathcal O(h^2)$ 
and $\mathcal O(h^{2.5})$, respectively.


\subsection{Example 2: A free final time FOCP with path and terminal constraints}

In this example, a free final time FOCP with path constraint is considered. 
Moreover, the state function is forced to lie on a circle at the final time:
\begin{equation}
\label{eq:ex2}
\begin{aligned}
\min J[{u}] &= \frac 1 2\int_0^{t_\text{f}\hspace{1pt}}\left[x^2(t)+u^2(t) \right]\text{\rm d}t\\
\text{s.t.}\ \ &{}_0^{\text{\tiny C}\hspace{1pt}}{\mathcal D}_t^\alpha  { x}(t) =-x(t)+u(t),\\
& {x}(0) =1, \\
& {u}(t) \ge 0.2, \ \ \ 0\le t \le t_\text{f}\hspace{1pt}, \\
& \left[x(t)-0.2\right]^2+\left[t-0.5\right]^2 \ge 0.25, \\
& \left[x(t_\text{f}\hspace{1pt})-0.2\right]^2+\left[t_\text{f}\hspace{1pt}-2\right]^2 = 0.04.
\end{aligned}
\end{equation}
Because of existence of path and terminal constraints, derivation of the optimality condition 
for this example is difficult. Consequently, solving this problem by indirect methods is cumbersome. 
We remind that in this paper we use a direct approach, which does not rely 
on the existence of optimality conditions. 

By applying the direct GL, TR and SI methods on this example, 
the obtained final time $t_\text{f}\hspace{1pt}$, the optimal objective 
value and CPU times for various values of $\alpha$ and $n$ are reported in Table~\ref{TAB2}.
It is worthwhile to note that problem \eqref{eq:ex2} is a free final time problem 
and has a nonsmooth solution. As a result, obtaining an accurate solution for this 
problem is not easily accessible. However, from Table~\ref{TAB2}, the precision 
and accuracy of the methods, especially the direct SI method, are concluded.
\begin{table}
\centering
\caption{(Example 2, with $\alpha=0.2, 0.4, \dots, 1.0$) CPU time, 
final time and value of the performance index, obtained through the direct GL, 
TR and SI methods with various values of $n$. All CPU times are in seconds.}
\label{TAB2}
\resizebox{\textwidth}{!}{
\begin{tabular}{lllllllllllllllllllll}
\hline
 &&  &  & \multicolumn{5}{c}{Direct GL method} &  & \multicolumn{5}{c}{Direct TR method} &  & 
\multicolumn{5}{c}{Direct SI method} \\
\cline{5-9}\cline{11-15}\cline{17-21}
\multicolumn{1}{c}{$\alpha$} && \multicolumn{1}{c}{ $n$} &&  \multicolumn{1}{c}{CPU} && 
\multicolumn{1}{c}{$t_\text{f}\hspace{1pt}$} &&  \multicolumn{1}{c}{$J_n$} &  &  \multicolumn{1}{c}{CPU} & & 
\multicolumn{1}{c}{$t_\text{f}\hspace{1pt}$} &  &\multicolumn{1}{c}{$J_n$} &  &  \multicolumn{1}{c}{CPU} & & 
\multicolumn{1}{c}{$t_\text{f}\hspace{1pt}$} & &\multicolumn{1}{c}{$J_n$} \\
\cline{1-1}\cline{3-3}\cline{5-5}\cline{7-7}\cline{9-9}\cline{11-11}\cline{13-13}\cline{15-15}\cline{17-17}
\cline{19-19}\cline{21-21}
0.2&&   31 &&0.2  &&1.859043&&0.309055&&0.1  &&1.859490&&0.318655&&0.1&&1.859530&&0.315002\\
0.2&&   61 &&0.3  &&1.859345&&0.309217&&0.2  &&1.859575&&0.313881&&0.3 &&1.859601&&0.312309\\
0.2&&   91 &&0.6  &&1.859440&&0.309338&&0.6  &&1.859595&&0.312419&&0.9 &&1.859614&&0.311426\\
0.2&&501&&73.1&&1.859599&&0.309759&&54.2 &&1.859628&&0.310313&&83.8&&1.859632&&0.310177\\
\hline
0.4&&   31 &&0.3  &&1.820755&&0.313741&&0.2  &&1.820827&&0.320906&&0.2 &&1.821028&&0.317973\\
0.4&&   61 &&0.4  &&1.820724&&0.314198&&0.4  &&1.820796&&0.3178  &&0.4 &&1.820789&&0.316745\\
0.4&&   91 &&0.7  &&1.820728&&0.314534&&0.9  &&1.820776&&0.316984&&0.8&&1.820761&&0.31639 \\
0.4&&501&&75.5&&1.820722&&0.315516&&85.3&&1.820731&&0.316007&&67.8&&1.820728&&0.315953\\ \hline
0.6&&   31 &&0.2  &&1.806128&&0.323821&&0.1  &&1.806192&&0.329454&&0.1  &&1.806075&&0.327235\\
0.6&&   61 &&0.3  &&1.806159&&0.324495&&0.3  &&1.805935&&0.327472&&0.4  &&1.805796&&0.326846\\
0.6&&   91 &&0.6  &&1.806034&&0.324961&&1.1  &&1.805920&&0.327014&&0.7 &&1.805890&&0.326683\\
0.6&&501&&91.1&&1.805853&&0.326182&&60.5&&1.805841&&0.326606&&77.4&&1.805833&&0.326589\\ \hline
0.8&&   31 &&0.1  &&1.801137&&0.335882&&0.1  &&1.801207&&0.339415&&0.1  &&1.801154&&0.337831\\
0.8&&   61 &&0.3  &&1.801029&&0.336271&&0.4  &&1.801109&&0.338177&&0.4  &&1.801077&&0.337766\\
0.8&&   91 &&0.5  &&1.801000&&0.336595&&0.5  &&1.801076&&0.337928&&0.9 	&&1.801053&&0.337733\\
0.8&&501&&63.2&&1.800979&&0.337453&&55.5&&1.801017&&0.337723&&  93.1&&1.801012&&0.337716\\ \hline
1.0&&   31 &&0.1  &&1.800442&&0.346938&&0.1  &&2.088769&&0.363078&&0.2   &&1.800840&&0.347474\\
1.0&&   61 &&0.3  &&1.800650&&0.34678 &&0.3  &&1.800884&&0.347631&&0.3  &&1.800904&&0.34732\\
1.0&&   91 &&0.5  &&1.800743&&0.346865&&0.6  &&1.800901&&0.347456&&0.8 &&1.800917&&0.347311\\
1.0&&501&&44.8&&1.800907&&0.347191&&71.8&&1.800939&&0.347304&&112.4&&1.800942&&0.347298\\ \hline
\end{tabular}
}
\end{table}
\begin{figure}
\includegraphics[width=\textwidth]{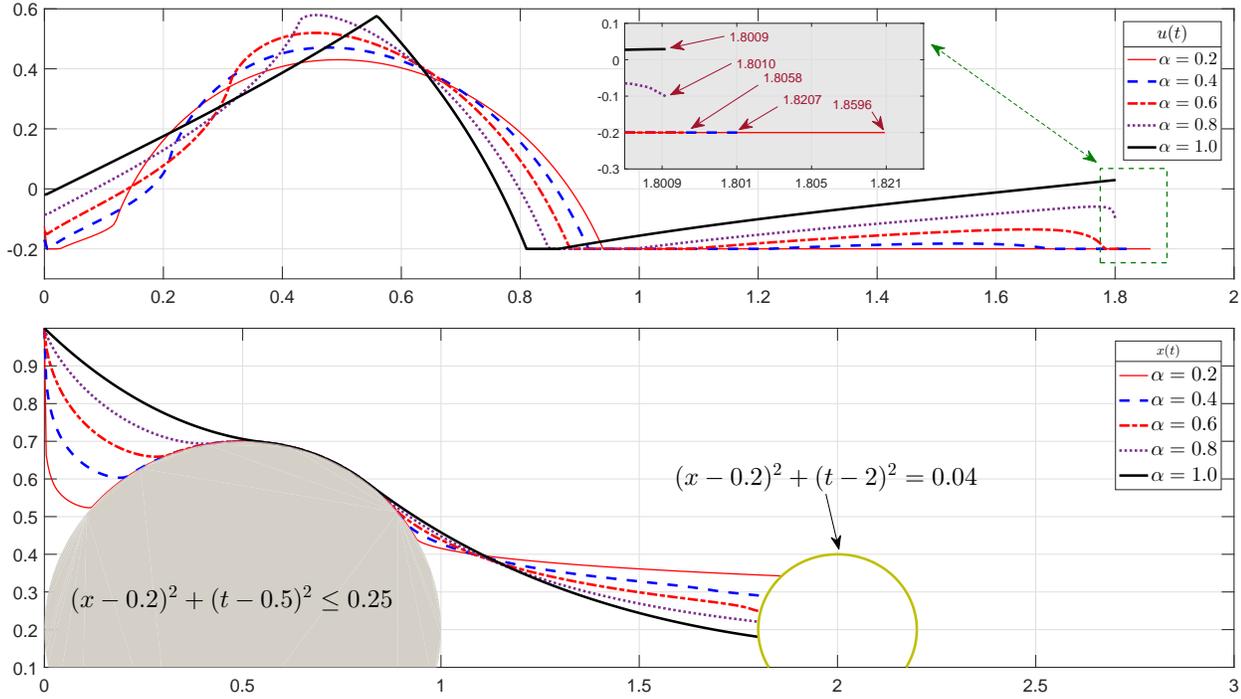}
\caption{(Example 2, with $\alpha=0.2,04,\dots,1.0$)  
The obtained control (above) and state (below) functions 
by the direct TR method.}
\end{figure}


\subsection{Example 3: A Bang-bang problem}

In this example, we apply our new numerical technique to the following 
bang-bang problem, which is treated in \cite{kamocki2015fractional}:
determine the state $x(t)$ and control $u(t)$ on the interval $t\in 
[0,2]$ that minimize the performance index
\[
J(x,u) = \int_0^{\mathop 2 } {\left[ {
x_{1}(t)- x_{2}(t)+ \mathop u(t) } \right]} \text{\rm d}t
\]
subject to the dynamic constraints
\begin{align*}
&{}_0^{\text{\tiny C}\hspace{1pt}}{\mathcal D}_t^\alpha x_{1}(t) =  x_{2}(t)- u(t),\\
&{}_0^{\text{\tiny C}\hspace{1pt}}{\mathcal D}_t^\alpha x_{2}(t) =  - u(t),
\end{align*}
and the initial conditions
\[
x_0=\left[ \begin{array}{ccc}
0 \\
1\end{array} \right].
\]
The exact solution to this problem is given by
\begin{equation*}
u_{*}(t) =
\begin{cases}
1  & \quad  t\in [0,1], \\
0  & \quad t\in [1,2],\\
\end{cases}
\qquad
x_{*}(t) =
\begin{bmatrix}
x_{*}^{1}(t) \\
x_{*}^{2}(t)
\end{bmatrix}
=
\left\{
\begin{array}{ll}
\begin{bmatrix}
-t\\
1-\frac{\sqrt{t}}{\Gamma(1.5)}
\end{bmatrix}, & \hbox{$t\in[0,1]$,} \\
\\
\begin{bmatrix}
\frac{\sqrt{t-1}}{\Gamma(1.5)}-1 \\
1-\frac{\sqrt{t}-\sqrt{t-1}}{\Gamma(1.5)}
\end{bmatrix}, & \hbox{$t\in[1,2]$.}
\end{array}
\right.
\end{equation*}
It is stressed that, in optimal control theory, the field 
of bang-bang optimal control problems is a classical topic. 
For bang-bang problems, the control function switches abruptly 
between its bounds. Bang-bang optimal control problems have received 
considerable attention, due to the difficulties arising in their numerical solutions.  
See \cite{Shamsi11OCAaM} and references therein. Here, we assess 
our methods with this bang-bang problem.

\begin{figure}
\includegraphics[width=\textwidth]{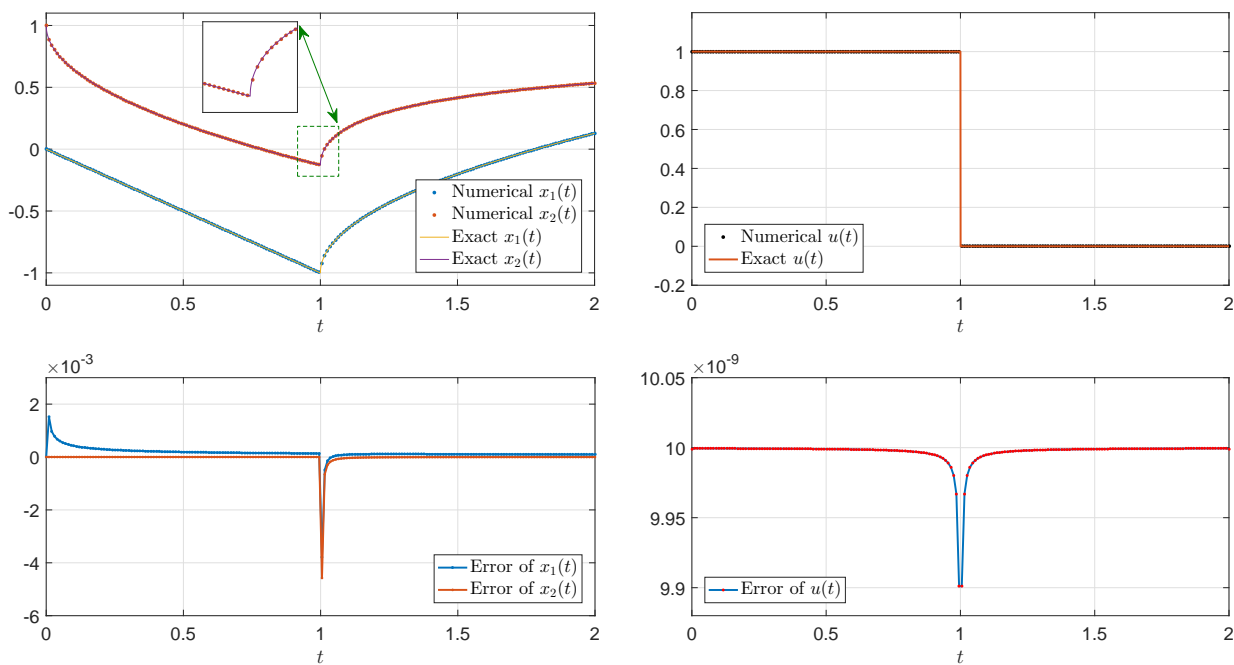}
\caption{(Example 3) Obtained solution via the direct TR method with $n=100$. 
Above Left: The exact and obtained states. Above Right: The exact and obtained control.  
Below Left: Error of the obtained states. Below Right: Error of the obtained control.}
\label{FIG4}
\end{figure}

\begin{figure}
\includegraphics[width=\textwidth]{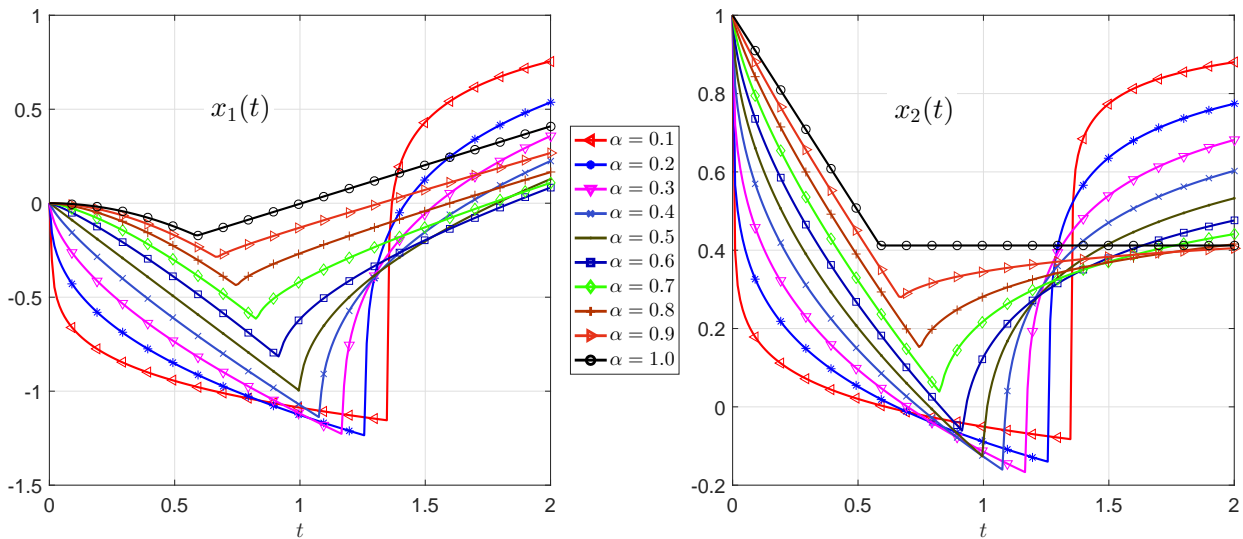}
\caption{(Example 3, with $\alpha=0.1,0.2,\dots,1.0$) 
The obtained state and control functions, 
by the direct TR method with $n=100$.}
\label{FIG5}
\end{figure}

The obtained control and state functions, by the direct TR method with $n=100$, 
are plotted in Figure~\ref{FIG4}. In addition, the errors of the obtained state 
and control are plotted in this figure too. It is seen that the control and state 
functions are accurately approximated by the direct TR method.
Moreover, we apply the direct TR method with $n=100$ on this example 
for $\alpha=0.1, 0.2, \dots, 0.9, 1.0$. The obtained state functions 
are plotted in Figure~\ref{FIG5}. We can see that in the all cases 
the control functions are bang-bang. In addition, the measured CPU times, 
optimal value of the objective function and switching times are reported 
in Table~\ref{TAB3}. To show the precision of the method, the results with 
$n=400$ are also reported in this table.
\begin{table}
\centering
\caption{(Example 3, with $\alpha=0.1,0.2,\dots,1.0$) 
The obtained CPU time, in seconds, value of the performance index and switching time, 
by the direct TR method with $n=100$ and $n=400$.}
\label{TAB3}
\begin{tabular}{lllllllll}
&  & \multicolumn{3}{c}{$n=100$} &  & \multicolumn{3}{c}{$n=400$} \\
\cline{3-5}\cline{7-9}
\multicolumn{1}{c}{$\alpha$} &  & \multicolumn{1}{c}{CPU} & \multicolumn{1}{c}{$J_n$} &
\multicolumn{1}{c}{Switch} &  & \multicolumn{1}{c}{CPU} & \multicolumn{1}{c}{$J_n$} &
\multicolumn{1}{c}{Switch} \\
\cline{1-1}\cline{3-5}\cline{7-9}
0.10&& 1.1&-0.14900  &1.34343   &&56.2&-0.14621  &1.34586   \\
0.20&& 0.5&-0.25034  &1.26262   &&31.3&-0.25109  &1.26065   \\
0.30&& 1.5&-0.32036  &1.17449   &&58.5&-0.32070  &1.17042   \\
0.40&& 0.5&-0.35859  &1.08080   &&36.6&-0.35912  &1.08521   \\
0.50&& 0.4&-0.37187  &1.00002   &&22.0&-0.37225  &1.00000   \\
0.60&& 0.5&-0.36618  &0.91910   &&42.6&-0.36644  &0.91478   \\
0.70&& 0.5&-0.34794  &0.83838   &&36.3&-0.34813  &0.83458   \\
0.80&& 1.4&-0.32337  &0.74496   &&34.1&-0.32343  &0.74937   \\
0.90&& 0.9&-0.29773  &0.67676   &&88.5&-0.29785  &0.66917   \\
1.00&& 1.4&-0.27611  &0.58389   &&51.2&-0.27613  &0.58395   \\ \hline
\end{tabular}	
\end{table}
We note that this problem is bang-bang and the control function 
is discontinuous and the state functions are nonsmooth. Naturally, 
the accuracy of a numerical method for such problems is lower in
comparison with smooth problems. Nevertheless, according to Table~\ref{TAB3}, 
the direct TR method provides solutions with reasonable accuracy 
for this bang-bang problem.


\subsection{Example 4: Optimal control of a fractional-order HIV-immune system with memory}

A HIV optimal control problem is now considered. 
The problem is stated in \cite{ding2012optimal} as follows:
\begin{eqnarray*}
&& {\rm min}\ J[u] = 500 [x_1^{2}(t_\text{f}\hspace{1pt}) 
+ x_3^{2}(t_\text{f}\hspace{1pt}) + x_4^{2}(t_\text{f}\hspace{1pt})] 
+ \int_0^{t_\text{f}\hspace{1pt} }\left(500 [x_1^{2}(\mathop t) + x_3^{2}(\mathop t) 
+ x_4^{2}(\mathop t)] + 0.005\,u_1^{2}(\mathop t)\right)\text{\rm d}t,\\
&&\qquad {}_0^{\text{\tiny C}\hspace{1pt}}{\mathcal D}_t^\alpha x_1(t)  
= - \mathop a\nolimits_1 \mathop x\nolimits_1(t)  - \mathop a\nolimits_2 
\mathop x\nolimits_1(t) \mathop x\nolimits_2(t) + \mathop a\nolimits_3 
\mathop a\nolimits_4 \mathop x\nolimits_4(t) (1 - \mathop u\nolimits_1(t) ),\\
&&\qquad {}_0^{\text{\tiny C}\hspace{1pt}}{\mathcal D}_t^\alpha x_2(t) 
= \frac{{\mathop a\nolimits_5 }}{{1 + \mathop x\nolimits_1(t) }} 
- \mathop a\nolimits_2 \mathop x\nolimits_1(t) \mathop x\nolimits_2(t) 
- \mathop a\nolimits_6 \mathop x\nolimits_2(t) + \mathop a\nolimits_7 
\left( {1 - \frac{{\mathop x\nolimits_2(t)  + \mathop x\nolimits_3(t)  
+ \mathop x\nolimits_4(t) }}{{\mathop a\nolimits_8 }}} \right)\mathop x\nolimits_2(t),\\
&&\qquad {}_0^{\text{\tiny C}\hspace{1pt}}{\mathcal D}_t^\alpha x_3(t) 
= \mathop a\nolimits_2 \mathop x\nolimits_1(t) \mathop x\nolimits_2(t) 
- \mathop a\nolimits_9 \mathop x\nolimits_3(t) - \mathop a\nolimits_6 \mathop x\nolimits_3(t),\\
&&\qquad {}_0^{\text{\tiny C}\hspace{1pt}}{\mathcal D}_t^\alpha x_4(t) 
= \mathop a\nolimits_9 \mathop x\nolimits_3(t) - \mathop a\nolimits_4 \mathop x\nolimits_4(t)\\
&&\qquad  \begin{array}{ll}
{\mathop x\nolimits_1 \left( 0 \right) = 0.049,}&{\mathop x\nolimits_2 \left( 0 \right) = 904,} \\
{\mathop x\nolimits_3 \left( 0 \right) = 0.034,}&{\mathop x\nolimits_4 \left( 0 \right) = 0.0042,}
\end{array}\\
\end{eqnarray*}
where the values of parameters $a_i$, $i=1,\dots,9$, are as shown in Table~\ref{TAB4}.
\begin{table}[ht]
\centering
\caption{Parameter values used in the optimal control of the fractional HIV-immune system.}
\label{TAB4}
\begin{tabular}{ccccccccccccccccc} \hline
$a_1$ &  & $a_2$ &  & $a_3$ &  & $a_4$ &  & $a_5$ &  & $a_6$ &  & $a_7$ &  & $a_8$ &  & $a_9$ \\
\cline{1-1}\cline{3-3}\cline{5-5}\cline{7-7}\cline{9-9}\cline{11-11}\cline{13-13}\cline{15-15}\cline{17-17}
2.4 &  & 2.4e-5 &  &1200  &  & 0.24  &  & 10  &  & 0.02  &  & 0.03 &  & 1500 &  &3e-3  \\ \hline
\end{tabular}
\end{table}
In \cite{ding2012optimal}, the authors used an indirect method to solve this FOCP. 
There, necessary optimality conditions for the problem are firstly derived, 
then the problem is solved using an iterative algorithm. Here, we solve 
the problem by the proposed direct methods.
 
\begin{figure}
\includegraphics[width=\textwidth]{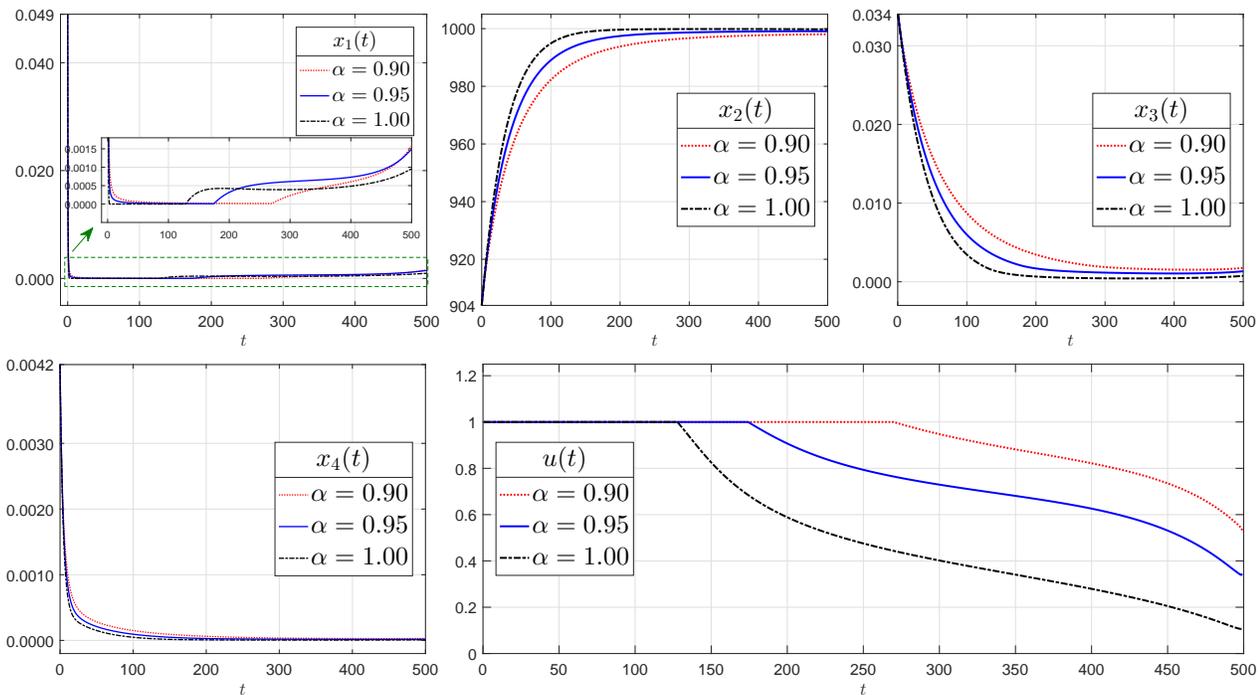}
\caption{(Example 4, with $\alpha=0.90,0.95,1.00$) 
Obtained states and control functions by the direct TR methods with $n=1000$.}
\label{FIG6}
\end{figure}

In Figure~\ref{FIG6}, we plot the obtained control and state functions, by applying 
the direct TR method with $n=1000$ to the problem for $\alpha=0.90,0.95,1.00$. 
Our results are in good agreement with those of \cite{ding2012optimal}. 
In addition, to report the precision of the methods, the values of the performance index 
obtained by the direct TR and SI methods are given in Table~\ref{TAB5}. 
Based on these results, we observe that, without deriving the optimality conditions, 
the problem can be solved with a reasonable accuracy.  
\begin{table}[ht]
\centering
\caption{(Example 4, with $\alpha=0.90,0.95,1.00$) Optimal value of the performance 
index obtained using the direct TR and SI methods with $n=500$ and $1000$.}
\label{TAB5}
\begin{tabular}{llllllllll} 
&  & \multicolumn{3}{c}{Direct TR method} &  &  & \multicolumn{3}{c}{Direct SI method} \\ 
\cline{3-5}\cline{8-10}
&  & $n=500$ &  & $n=1000$ &  &  & $n=500$ &  & $n=1000$ \\ 
\cline{1-1}\cline{3-3}\cline{5-5}\cline{8-8}\cline{10-10}
$\alpha=0.90$ &  & 22.70 &  & 22.65 &  &  & 22.63 &  & 22.61 \\ 
$\alpha=0.90$ &  & 18.07 &  & 18.03 &  &  & 18.01 &  & 17.99 \\ 
$\alpha=1.00$ &  & 14.72 &  & 14.39 &  &  & 14.32 &  & 14.31 \\ \hline
\end{tabular}	
\end{table}


\section{Conclusion}
\label{sec:6}

In this paper, three direct methods based on Gr\"{u}nwald-Letnikov, trapezoidal 
and Simpson approximation formulas are presented for the numerical solution 
of a general class of fractional optimal control problems. 
Optimality conditions are not needed in these methods. Thus, they can 
be applied to any type of fractional optimal control problem.

The proposed direct methods are illustrated in four academic and practical test problems 
and the results confirm that our methods are reliable.
According to our numerical results, for problems with a smooth solution, 
the accuracy of the direct Simpson method is superior and for problems 
with discontinuous or a nonsmooth solution, the accuracy is satisfactory. 
Moreover, by providing the gradient of the objective function and the 
Jacobian of the constraints, the CPU time is significantly reduced. 
We conclude that the direct methods here proposed are simple, 
reliable, reasonably accurate and fast for solving fractional optimal control problems.

As further research works, we can refer to costate estimation and combining 
mesh generation techniques with the presented method to improve the accuracy 
and speed of the methods in solving problems with nonsmooth solutions. 


\section*{Acknowledgements}

Torres has been partially supported by FCT through the R\&D Unit CIDMA (UID/MAT/04106/2013) 
and TOCCATA project PTDC/EEI-AUT/2933/2014 funded by FEDER and COMPETE 2020.



\end{document}